\documentstyle[12pt,twoside,leqno]{amsart}

\newtheorem{Theorem}{Theorem}[section]
\newtheorem{Proposition}[Theorem]{Proposition}
\newtheorem{Lemma}[Theorem]{Lemma}
\newtheorem{Corollary}[Theorem]{Corollary}

\theoremstyle{remark}
\newtheorem{Remark}[Theorem]{Remark}

\theoremstyle{definition}
\newtheorem{Definition}[Theorem]{Definition}

\newtheorem{Counterexample}[Theorem]{Counterexample}

\newtheorem{Question}[Theorem]{Question}

\begin{document}
\bibliographystyle{plain}
\title{Newtonian Lorentz metric spaces}
\author[\c{S}. Costea]{\c{S}erban Costea}
\address{\c{S}. Costea\\
EPFL SB MATHGEOM\\
Section de Math\'{e}matiques\\
Station 8\\
CH-1015 Lausanne, Switzerland}
\email{serban.costea@@epfl.ch}
\author[M. Miranda Jr.]{Michele Miranda Jr.}
\address{M. Miranda Jr.\\
Department of Mathematics\\
University of Ferrara\\
Via Machiavelli 35\\
I-44121 Ferrara, Italy}
\email{michele.miranda@@unife.it}


\begin{abstract}

 This paper studies Newtonian Sobolev-Lorentz spaces. We prove that these spaces are Banach. We also study the global $p,q$-capacity and the $p,q$-modulus of families of rectifiable curves. Under some additional assumptions (that is, $X$ carries a doubling measure and a weak Poincar\'{e} inequality), we show that when $1 \le q<p$ the Lipschitz functions are dense in these spaces; moreover, in the same setting we show that the $p,q$-capacity is Choquet provided that $q>1.$ We also provide a counterexample to the density result in the Euclidean setting when $1 < p \le n$ and $q=\infty.$
\end{abstract}

\maketitle
\section{Introduction}
 In this paper $(X,d)$ is a complete metric space endowed with a nontrivial Borel regular measure $\mu.$ We assume that $\mu$ is finite and nonzero on nonempty bounded open sets. In particular, this implies that the measure $\mu$ is $\sigma$-finite. Further restrictions on the space $X$ and on the measure $\mu$ will be imposed later.

 The Sobolev-Lorentz relative $p,q$-capacity was studied in the Euclidean setting by  Costea \cite{Cos} and  Costea-Maz'ya \cite{CosMaz}.
 The Sobolev $p$-capacity was studied by Maz'ya \cite{Maz} and Heinonen-Kilpel\"{a}inen-Martio
 \cite{HKM} in ${\mathbf{R}}^n$
 and by Costea \cite{Cos1} and Kinnunen-Martio \cite{KiM1} and \cite{KiM2} in metric spaces. The relative Sobolev $p$-capacity in metric spaces was introduced by J. Bj{\"o}rn in \cite{Bjo} when studying the boundary continuity properties of quasiminimizers.

After recalling the definition of $p,q$-Lorentz spaces in metric spaces, we
study some useful property of the $p,q$-modulus of families of curves
needed to give the notion of $p,q$-weak upper gradients.  Then,
following the approach of Shanmugalingam in \cite{Sha1} and \cite{Sha2}, we generalize
the notion of Newtonian Sobolev spaces to the Lorentz setting.
 There are several other definitions of Sobolev-type spaces in the metric setting
 when $p=q$; see Haj{\l}asz \cite{Haj1}, Heinonen-Koskela \cite{HeK}, Cheeger \cite{Che}, and Franchi-Haj{\l}asz-Koskela \cite{FHK}. It has been shown that under reasonable hypotheses, the majority of these definitions yields the same space; see Franchi-Haj{\l}asz-Koskela \cite{FHK} and Shanmugalingam \cite{Sha1}.

We prove that these spaces are Banach. In order to this, we develop a theory of
Sobolev $p,q$-capacity. Some of the ideas used here when proving the properties of the $p,q$-capacity follow Kinnunen-Martio \cite{KiM1} and \cite{KiM2} and Costea \cite{Cos1}. We also use this theory to prove that, in the case $1 \le q<p,$
Lipschitz functions are dense in the Newtonian Sobolev-Lorentz space if the space $X$ carries a doubling measure $\mu$ and a weak $(1, L^{p,q})$-Poincar\'{e} inequality. Newtonian Banach-valued Sobolev-Lorentz spaces were studied by Podbrdsky in \cite{Pod}.

 We prove that under certain restrictions (when $1< q \le p$ and the space $(X, d)$ carries a doubling measure $\mu$ and a certain weak Poincar\'{e} inequality) this capacity is a Choquet set function.

 We recall the standard notation and definitions to be used throughout
 this paper. We denote by $B(x,r)=\{ y \in X: d(x,y)<r \}$  the open ball with center $x \in X$ and radius $r>0,$ while $\overline{B}(x,r)=\{ y \in X: d(x,y) \le r \}$
 is the closed ball with center $x \in X$ and radius $r>0.$ For a positive
 number $\lambda,$ $\lambda B(a,r)=B(a, \lambda r)$ and $\lambda \overline{B}(a,r)=\overline{B}(a, \lambda r).$

 Throughout this paper, $C$ will denote a positive constant whose value
 is not necessarily the same at each occurrence; it may vary even within
 a line. $C(a,b, \ldots)$ is a constant that depends only on the parameters
 $a,b, \ldots.$
For $E \subset X,$ the boundary, the closure, and the complement of $E$ with respect
 to $X$ will be denoted by $\partial E,$ $\overline{E},$ and $X \setminus E,$
 respectively; $\mbox{diam }E$ is the diameter of $E$ with
 respect to the metric $d$.

 \section{Lorentz spaces}

 Let $f: X \rightarrow [-\infty, \infty]$ be a $\mu$-measurable function.
 We define $\mu_{[f]},$ the \textit{distribution function}
 of $f$ as follows (see Bennett-Sharpley \cite[Definition II.1.1]{BS}):
 $$\mu_{[f]}(t)=\mu(\{x \in X: |f(x)| > t \}), \qquad t \ge 0.$$
 We define $f^{*},$  the \textit{nonincreasing rearrangement} of $f$ by
$$f^{*}(t)=\inf\{v: \mu_{[f]}(v) \le t \}, \quad t \ge 0.$$
(See Bennett-Sharpley \cite[Definition II.1.5]{BS}.)
 We note that $f$ and $f^{*}$ have the same distribution function.
 For every positive $\alpha$ we have $$(|f|^{\alpha})^{*}=(|f|^{*})^{\alpha}$$
 and if $|g|\le |f|$ $\mu$-almost everywhere on $X,$ then $g^{*}\le f^{*}.$
 (See \cite[Proposition II.1.7]{BS}.)
 We also define $f^{**}$, the \textit{maximal function} of $f^{*}$ by
  $$f^{**}(t)=m_{f^{*}}(t)=\frac{1}{t} \int_{0}^{t} f^{*}(s) ds, \quad t >0.$$
 (See \cite[Definition II.3.1]{BS}.)

Throughout the paper, we denote by $p'$ the H\"{o}lder
conjugate of $p \in [1,\infty]$.

\smallskip

The \textit{Lorentz space} $L^{p,q}(X, \mu),$
$1<p<\infty,$ $1\le q\le \infty,$ is defined as follows:
$$L^{p,q}(X, \mu)= \{f: X \rightarrow [-\infty, \infty]: f \mbox { is $\mu$-measurable, }
||f||_{L^{p,q}(X, \mu)}<\infty \},$$
where
$$||f||_{L^{p,q}(X,\mu)}=||f||_{p,q}=\left\{ \begin{array}{lc}
\left( \displaystyle{\int_{0}^{\infty} (t^{1/p}f^{*}(t))^q \, \frac{dt}{t}}
\right)^{1/q}, & 1 \le q < \infty, \\
\\
\sup\limits_{t>0} t \mu_{[f]}(t)^{1/p}=\sup\limits_{s>0}
s^{1/p} f^{*}(s), & q=\infty.
\end{array}
\right.
$$
(See Bennett-Sharpley \cite[Definition IV.4.1]{BS} and Stein-Weiss
\cite[p.\ 191]{SW}.)

\smallskip

If $1 \le q\le p,$ then $||\cdot||_{L^{p,q}(X, \mu)}$
represents a norm, but for $p < q \le \infty$ it represents a
quasinorm, equivalent to the norm $||\cdot||_{L^{(p,q)}(X, \mu)},$ where
$$||f||_{L^{(p,q)}(X, \mu)}=||f||_{(p,q)}=\left\{ \begin{array}{lc}
\left(\displaystyle{ \int_{0}^{\infty} (t^{1/p}f^{**}(t))^q \, \frac{dt}{t} }\right)^{1/q}, & 1 \le q < \infty, \\
\\
\sup\limits_{t>0} t^{1/p} f^{**}(t), & q=\infty.
\end{array}
\right.
$$
(See \cite[Definition IV.4.4]{BS}.) Namely, from \cite[Lemma IV.4.5]{BS} we have that
$$||f||_{L^{p,q}(X,\mu)} \le ||f||_{L^{(p,q)}(X,\mu)} \le p' ||f||_{L^{p,q}(X,\mu)}$$
for every $q \in [1, \infty]$ and every $\mu$-measurable function
$f:X \rightarrow [-\infty, \infty].$

It is known that $(L^{p,q}(X,\mu),
||\cdot||_{L^{p,q}(X,\mu)})$ is a Banach space for
$1\le q \le p,$ while $(L^{p,q}(X,\mu),
||\cdot||_{L^{(p,q)}(X,\mu)})$ is a Banach space
for $1<p< \infty,$ $1\le q \le \infty.$ In addition, if the measure $\mu$ is nonatomic, the aforementioned Banach spaces are reflexive when $1<q<\infty.$ (See Hunt \cite[p.\ 259-262]{Hu} and Bennett-Sharpley \cite[Theorem IV.4.7 and Corollaries I.4.3 and IV.4.8]{BS}.) (A measure $\mu$ is called \textit{nonatomic} if for every measurable set $A$ of positive measure there exists a measurable set $B \subset A$ such that $0<\mu(B)<\mu(A)$.)

\smallskip

\begin{Definition} \label{definition of abs continuous norm} (See \cite[Definition I.3.1]{BS}.) Let $1 < p < \infty$ and $1 \le q \le \infty.$ Let $Y = L^{p,q}(X, \mu).$ A function $f$ in $Y$ is said to have absolutely continuous norm in $Y$ if and only if $||f \chi_{E_k}||_{Y} \rightarrow 0$ for every sequence $E_k$ of $\mu$-measurable sets satisfying $E_{k} \rightarrow \emptyset$ $\mu$-almost everywhere.
\end{Definition}

Let $Y_{a}$ be the subspace of $Y$ consisting of functions of absolutely continuous
norm and let $Y_{b}$ be the closure in $Y$ of the set of simple functions. It is known
that $Y_{a} = Y_{b}$ whenever $1 \le q \le \infty.$  (See Bennett-Sharpley \cite[Theorem I.3.13]{BS}.) Moreover, since $(X, \mu)$ is a $\sigma$-finite measure space, we have $Y_{b} = Y$ whenever $1\le q < \infty.$ (See Hunt \cite[p.\ 258-259]{Hu}.)

We recall (see Costea \cite{Cos}) that in the Euclidean setting (that is, when $\mu=m_n$ is the $n$-dimensional Lebesgue measure and $d$ is the Euclidean distance on ${\Bbb{R}}^n$) we have $Y_a \neq Y$ for $Y=L^{p,\infty}(X, m_n)$ whenever $X$ is an open subset of ${\Bbb{R}}^n.$ Let $X=B(0,2) \setminus \{0\}.$ As in Costea \cite{Cos} we define $u: X \rightarrow {\Bbb{R}},$
\begin{eqnarray}\label{counterexample wk Lp is not abs cont}
 u(x)=\left\{
\begin{array}{cl}
 |x|^{-\frac{n}{p}} & \mbox{if $0<|x|<1$}\\
 0 & \mbox{if $1 \le |x| \le 2.$}
 \end{array}
\right.
\end{eqnarray}
It is easy to see that $u \in L^{p,\infty}(X, m_n)$ and moreover,
\begin{equation*}
||u\chi_{B(0,\alpha)}||_{L^{p,\infty}(X, m_n)}=||u||_{L^{p,\infty}(X, m_n)}=m_n(B(0,1))^{1/p}
\end{equation*}
for every $\alpha>0.$ This shows that $u$ does not have absolutely
continuous weak $L^{p}$-norm and therefore $L^{p,\infty}(X, m_n)$
does not have absolutely continuous norm.

\begin{Remark} \label{relation between Lpr and Lps} It is also known (see \cite[Proposition IV.4.2]{BS}) that for every $p \in (1,\infty)$ and
$1\le r<s\le \infty$ there exists a constant $C(p,r,s)$ such that
\begin{equation}\label{relation between Lpr and Lps norm}
||f||_{L^{p,s}(X, \mu)} \le C(p,r,s)
||f||_{L^{p,r}(X, \mu)}
\end{equation}
for all measurable functions $f \in L^{p,r}(X, \mu).$
In particular, the embedding $L^{p,r}(X, \mu) \hookrightarrow
L^{p,s}(X, \mu)$ holds.
\end{Remark}

\begin{Remark} \label{identities with powers of norms} Via Bennett-Sharpley \cite[Proposition II.1.7 and Definition IV.4.1]{BS} it is easy to see
that for every $p \in (1, \infty),$ $q \in [1, \infty]$ and $0< \alpha \le \min(p,q),$
we have
$$||f||_{L^{p,q}(X, \mu)}^{\alpha}=||f^{\alpha}||_{L^{\frac{p}{\alpha},\frac{q}{\alpha}}(X, \mu)}$$
for every nonnegative function $f \in L^{p,q}(X, \mu).$
\end{Remark}

\subsection{\bf{The subadditivity and superadditivity of the Lorentz quasinorms}}

We recall the known results and present new results concerning the superadditivity and the subadditivity of the Lorentz $p,q$-quasinorm. For the convenience of the reader, we will provide proofs for the new results and for some of the known results.

The superadditivity of the Lorentz $p,q$-norm in the case $1 \le q \le p$ was stated in  Chung-Hunt-Kurtz \cite[Lemma 2.5]{CHK}.

\begin{Proposition} \label{superadd q le p} (See \cite[Lemma 2.5]{CHK}.)
Let $(X, \mu)$ be a measure space. Suppose that $1 \le q \le p.$ Let $\{E_i\}_{i \ge 1}$ be
a collection of pairwise disjoint $\mu$-measurable subsets of $X$ with $E_{0}=\cup_{i \ge 1} E_i$ and let $f \in L^{p,q}(X, \mu).$ Then
\begin{equation*}
 \sum_{i \ge 1} || \chi_{E_i} f||_{L^{p,q}(X, \mu)}^p \le || \chi_{E_0} f||_{L^{p,q}(X, \mu)}^p.
\end{equation*}
\end{Proposition}

A similar result concerning the superadditivity was obtained in Costea-Maz'ya \cite[Proposition 2.4]{CosMaz} for the case $1<p < q<\infty$ when $X=\Omega$ was an open set in ${\Bbb{R}}^n$ and $\mu$ was an arbitrary measure. That result is valid for a general measure space $(X, \mu).$

\begin{Proposition} \label{superadd p le q} Let $(X, \mu)$ be a measure space. Suppose that $1<p < q<\infty.$ Let $\{E_i\}_{i \ge 1}$ be
a collection of pairwise disjoint $\mu$-measurable subsets of $X$ with $E_{0}=\cup_{i \ge 1} E_i$ and let $f \in L^{p,q}(X, \mu).$ Then
\begin{equation*}
 \sum_{i \ge 1} || \chi_{E_i} f||_{L^{p,q}(X, \mu)}^q \le || \chi_{E_0} f||_{L^{p,q}(X, \mu)}^q.
\end{equation*}
\end{Proposition}

\begin{pf} We mimic the proof of Proposition 2.4 from Costea-Maz'ya \cite{CosMaz}. We replace $\Omega$ with $X.$
\end{pf}

We have a similar result for the subadditivity of the Lorentz $p,q$-quasinorm.
When $1 < p < q \le \infty$ we obtain a result that generalizes Theorem 2.5 from Costea \cite{Cos}.

\begin{Proposition} \label{subadd p le q} Let $(X, \mu)$ be a measure space. Suppose that
$1 < p < q \le \infty.$ Suppose $f_i, i=1,2, \ldots$ is a sequence of functions in $L^{p,q}(X, \mu)$ and let $f_{0}=\sup_{i \ge 1} |f_i|.$  Then

 $$||f_{0}||_{L^{p,q}(X, \mu)}^p \le \sum_{i=1}^{\infty} ||f_i||_{L^{p,q}(X, \mu)}^p.$$

\end{Proposition}

\begin{pf} Without loss of generality we can assume that all the functions $f_i, i=1,2, \ldots$ are nonnegative. We have to consider two cases, depending on whether $p<q<\infty$ or $q=\infty.$

Let $\mu_{[f_i]}$ be the distribution function of $f_i$ for
$i=0,1,2, \ldots.$ It is easy to see that
\begin{equation} \label{the distrib fn of f0 less than the sum of distrib fn fi}
\mu_{[f_0]}(s) \le \sum_{i=1}^{\infty} \mu_{[f_i]}(s)
\mbox{ for every } s\ge 0.
\end{equation}

Suppose that $p<q<\infty.$ We have (see Kauhanen-Koskela-Mal{\'{y}} \cite[Proposition 2.1]{KKM})
\begin{equation} \label{KKM formula q finite}
||f_i||_{L^{p,q}(X, \mu)}^{p}=\left(p \int_{0}^{\infty} s^{q-1}
{\mu_{[f_i]}(s)}^{\frac{q}{p}} ds\right)^{\frac{p}{q}}
\end{equation}
for $i=0,1,2, \ldots.$ From this and
(\ref{the distrib fn of f0 less than the sum of distrib fn fi}) we obtain
\begin{eqnarray*}
||f_0||_{L^{p,q}(\Omega, \mu)}^{p} & = & \left(p \int_{0}^{\infty}
s^{q-1}  \mu_{[f_0]}(s)^{\frac{q}{p}} ds\right)^{\frac{p}{q}} \le \sum_{i \ge 1} \left(p \int_{0}^{\infty} s^{q-1}
{\mu_{[f_i]}(s)}^{\frac{q}{p}} ds\right)^{\frac{p}{q}} \\
& = &\sum_{i \ge 1} ||f_i||_{L^{p,q}(\Omega, \mu)}^{p}.
\end{eqnarray*}

Now, suppose that $q=\infty.$ From
(\ref{the distrib fn of f0 less than the sum of distrib fn fi}) we obtain
\begin{equation*}
s^p \, \mu_{[f_0]}(s) \le \sum_{i \ge 1} (s^p \, \mu_{[f_i]}(s)) \mbox{ for every $s>0,$}
\end{equation*}
which implies
\begin{equation}\label{p, inf equality 1}
s^p \, \mu_{[f_0]}(s) \le \sum_{i \ge 1} ||f_i||_{L^{p, \infty}(X, \mu)}^p \mbox{ for every $s>0.$}
\end{equation}
By taking the supremum over all $s>0$ in (\ref{p, inf equality 1}), we get the desired conclusion. This finishes the proof.

\end{pf}

We recall a few results concerning Lorentz spaces.

\begin{Theorem} \label{Quasinorm p<q leq infty power property} (See \cite[Theorem 2.6]{Cos}.)
Suppose $1<p<q \le \infty$ and $\varepsilon \in (0,1).$ Let $f_1, f_2 \in L^{p,q}(X, \mu).$ We
denote $f_3=f_1+f_2.$ Then $f_3 \in L^{p,q}(X, \mu)$ and
$$||f_3||_{L^{p,q}(X, \mu)}^p \le
(1-\varepsilon)^{-p}||f_1||_{L^{p,q}(X, \mu)}^p+ \varepsilon^{-p}
||f_2||_{L^{p,q}(X, \mu)}^p.$$

\end{Theorem}

\begin{pf} The proof of Theorem 2.6 from Costea \cite{Cos} carries verbatim. We replace $\Omega$ with $X.$

\end{pf}

Theorem \ref{Quasinorm p<q leq infty power property} has an useful corollary.

\begin{Corollary} \label{Convergence of the quasinorm of the sequence} (See \cite[Corollary 2.7]{Cos}.)
Suppose $1<p<\infty$ and $1\le q \le \infty.$ Let $f_k$ be a sequence of functions in
$L^{p,q}(X, \mu)$ converging to $f$ with respect to
the $p,q$-quasinorm and pointwise $\mu$-almost everywhere in $X.$ Then
$$\lim_{k \rightarrow
\infty} ||f_k||_{L^{p,q}(X, \mu)}=||f||_{L^{p,q}(X, \mu)}.$$
\end{Corollary}

\begin{pf} The proof of Corollary 2.7 from Costea \cite{Cos} carries verbatim. We replace $\Omega$ with $X.$

\end{pf}

 \section{p,q-modulus of the path family}

 In this section, we establish some results about $p,q$-modulus of families of curves. Here $(X, d, \mu)$ is a metric measure space.
 We say that a curve $\gamma$ in $X$ is rectifiable if it has finite length. Whenever $\gamma$ is rectifiable, we use the arc length parametrization $\gamma: [0, \ell(\gamma)] \rightarrow X,$ where $\ell(\gamma)$ is the length of the curve $\gamma.$

 Let $\Gamma_{\rm{rect}}$ denote the family of all nonconstant rectifiable curves in $X.$ It may well be that $\Gamma_{\rm{rect}}=\emptyset$, but we will be interested
 in metric spaces for which $\Gamma_{\rm{rect}}$ is sufficiently large.

 \begin{Definition} For $\Gamma \subset \Gamma_{\rm{rect}},$ let $F(\Gamma)$ be the family of all Borel measurable functions $\rho: X \rightarrow [0, \infty]$ such that
 $$\int_{\gamma} \rho \ge 1 \mbox{ for every } \gamma \in \Gamma.$$

 Now for each $1<p<\infty$ and $1 \le q \le \infty$ we define
 $${\mathrm{Mod}}_{p,q}(\Gamma)= \inf_{\rho \in F(\Gamma)} ||\rho||_{L^{p,q}(X,\mu)}^p.  $$

 \end{Definition}

 The number ${\mathrm{Mod}}_{p,q}(\Gamma)$ is called the \textit{p,q-modulus} of the family $\Gamma.$

 \subsection{Basic properties of the \textit{p,q}-modulus} Usually, a modulus is a monotone and subadditive set function. The following theorem will show, among other things, that this is true in the case of the $p,q$-modulus.

\begin{Theorem} \label{basic properties of pq modulus} Suppose $1<p<\infty$ and $1 \le q \le \infty.$
The set function $\Gamma \rightarrow {\mathrm{Mod}}_{p,q}(\Gamma),$ $\Gamma \subset \Gamma_{\rm{rect}},$ enjoys the following properties:

\par{\rm{(i)}} ${\mathrm{Mod}}_{p,q}(\emptyset)=0.$

\par {\rm{(ii)}}  If $\Gamma_{1} \subset \Gamma_{2},$ then ${\mathrm{Mod}}_{p,q}(\Gamma_{1}) \le {\mathrm{Mod}}_{p,q}(\Gamma_{2}).$

\par {\rm{(iii)}} Suppose $1 \le q \le p.$ Then $${\mathrm{Mod}}_{p,q}(\bigcup_{i=1}^{\infty} \Gamma_{i})^{q/p} \le \sum_{i=1}^{\infty} {\mathrm{Mod}}_{p,q}(\Gamma_{i})^{q/p}.$$

\par {\rm{(iv)}} Suppose $p < q \le \infty.$ Then $${\mathrm{Mod}}_{p,q}(\bigcup_{i=1}^{\infty} \Gamma_{i}) \le \sum_{i=1}^{\infty} {\mathrm{Mod}}_{p,q}(\Gamma_{i}).$$

\end{Theorem}

\begin{pf} (i) ${\mathrm{Mod}}_{p,q}(\emptyset)=0$ because $\rho \equiv 0 \in F(\emptyset).$

\vskip 2mm

 \par (ii) If $\Gamma_{1} \subset \Gamma_{2},$ then $F(\Gamma_2) \subset F(\Gamma_1)$ and hence ${\mathrm{Mod}}_{p,q}(\Gamma_{1}) \le {\mathrm{Mod}}_{p,q}(\Gamma_{2}).$

\vskip 2mm

 \par (iii) Suppose that $1 \le q \le p.$ The case $p=q$ corresponds to the $p$-modulus
 and the claim certainly holds in that case. (See for instance Haj{\l}asz
 \cite[Theorem 5.2 (3)]{Haj2}.) So we can look at the case $1\le q<p.$

  We can assume without loss of generality that $$\sum_{i=1}^{\infty} {\mathrm{Mod}}_{p,q}(\Gamma_{i})^{q/p}<\infty.$$

 Let $\varepsilon>0$ be fixed. Take $\rho_i \in F(\Gamma_i)$ such that $$||\rho_i||_{L^{p,q}(X,\mu)}^q <{\mathrm{Mod}}_{p,q}(\Gamma_{i})^{q/p}+\varepsilon 2^{-i}.$$ Let $\rho:=(\sum_{i=1}^\infty \rho_i^q)^{1/q}.$ We notice via Bennett-Sharpley \cite[Proposition II.1.7 and Definition IV.4.1]{BS} and Remark \ref{identities with powers of norms} applied with $\alpha=q$ that
  \begin{equation} \label{relations between quasinorms to different powers}
  \rho_i^q \in L^{\frac{p}{q}, 1}(X, \mu) \mbox{ and } ||\rho_i^q||_{L^{\frac{p}{q}, 1}(X, \mu)}=||\rho_i||_{L^{p,q}(X, \mu)}^q.
 \end{equation}
 for every $i=1,2, \ldots.$
  By using (\ref{relations between quasinorms to different powers}) and Remark \ref{identities with powers of norms} together with the definition of $\rho$ and the fact that $||\cdot||_{L^{{\frac{p}{q},1}}(X,\mu)}$ is a norm when $1 \le q \le p,$ it follows that $\rho \in F(\Gamma)$ and
 $${\mathrm{Mod}}_{p,q}(\Gamma_{i})^{q/p} \le ||\rho||_{L^{p,q}(X, \mu)}^q \le \sum_{i=1}^{\infty} ||\rho_i||_{L^{p,q}(X, \mu)}^q <\sum_{i=1}^{\infty} {\mathrm{Mod}}_{p,q}(\Gamma_{i})^{q/p}+2 \varepsilon.$$ Letting $\varepsilon \rightarrow 0,$ we complete the proof when $1 \le q \le p.$

 \vskip 2mm

 \par (iv) Suppose now that $p<q \le \infty.$ We can assume without loss of generality that
 $$\sum_{i=1}^{\infty} {\mathrm{Mod}}_{p,q}(\Gamma_{i})<\infty.$$
 Let $\varepsilon>0$ be fixed. Take $\rho_i \in F(\Gamma_i)$ such that $$||\rho_i||_{L^{p,q}(X,\mu)}^p <{\mathrm{Mod}}_{p,q}(\Gamma_{i})+\varepsilon 2^{-i}.$$ Let $\rho:=\sup_{i \ge 1} \rho_i.$ Then $\rho \in F(\Gamma).$ Moreover, from Proposition \ref{subadd p le q} it follows that $\rho \in L^{p,q}(X, \mu)$ and
 $${\mathrm{Mod}}_{p,q}(\Gamma) \le ||\rho||_{L^{p,q}(X, \mu)}^p \le \sum_{i=1}^{\infty} ||\rho_i||_{L^{p,q}(X, \mu)}^p <\sum_{i=1}^{\infty} {\mathrm{Mod}}_{p,q}(\Gamma_{i})+2 \varepsilon.$$ Letting $\varepsilon \rightarrow 0,$ we complete the proof when $p < q \le \infty.$

\end{pf}

 So we proved that the modulus is a monotone function. Also, the shorter the curves, the larger the modulus. More precisely, we have:

 \begin{Lemma} Let $\Gamma_{1}, \Gamma_2 \subset \Gamma_{\rm{rect}}.$ If each curve
 $\gamma \in \Gamma_1$ contains a subcurve that belongs to $\Gamma_2,$ then
 ${\mathrm{Mod}}_{p,q}(\Gamma_{1})\le{\mathrm{Mod}}_{p,q}(\Gamma_{2}).$
 \end{Lemma}

 \begin{pf} $F(\Gamma_2) \le F(\Gamma_1).$

 \end{pf}

 The following theorem provides an useful characterization of path families
 that have $p,q$-modulus zero.

 \begin{Theorem} \label{char of families of curves of pq modulus zero} Let $\Gamma \subset \Gamma_{\rm{rect}}.$ Then ${\mathrm{Mod}}_{p,q}(\Gamma)=0$ if and only if there exists a Borel measurable function $0 \le \rho \in L^{p,q}(X, \mu)$ such that $\int_{\gamma} \rho=\infty$ for every $\gamma \in \Gamma.$

 \end{Theorem}

 \begin{pf} Sufficiency. We notice that $\rho/n \in F(\Gamma)$ for every $n$ and hence $${\mathrm{Mod}}_{p,q}(\Gamma) \le \lim_{n \rightarrow \infty} ||\rho/n||_{L^{p,q}(X,\mu)}^p=0.$$

  Necessity. There exists $\rho_i \in F(\Gamma)$ such that $||\rho_i||_{L^{(p,q)}(X,\mu)}<2^{-i}$ and $\int_{\gamma} \rho_{i} \ge 1$ for every $\gamma \in \Gamma.$ Then $\rho:=\sum_{i=1}^{\infty} \rho_i$ has the required properties.

 \end{pf}

 \begin{Corollary} Suppose $1<p<\infty$ and $1 \le q \le \infty$ are given. If $0 \le g \in L^{p,q}(X, \mu)$ is Borel measurable, then $\int_{\gamma} g <\infty$ for $p,q$-almost every $\gamma \in \Gamma_{\rm{rect}}.$

 \end{Corollary}

 The following theorem is also important.

 \begin{Theorem} \label{uk conv to u in Lpq implies uk conv to u on curves} Let $u_k: X \rightarrow \overline{\Bbb{R}}=[-\infty, \infty]$ be a sequence of Borel functions which converge to a Borel function
 $u: X \rightarrow \overline{\Bbb{R}}$ in $L^{p,q}(X,\mu).$ Then there is a subsequence $(u_{k_j})_j$ such that
 $$\int_{\gamma} |u_{k_j}-u| \rightarrow 0 \mbox { as } j \rightarrow \infty,$$
 for $p,q$-almost every curve $\gamma \in \Gamma_{\rm{rect}}.$

 \end{Theorem}

 \begin{pf} We follow Haj{\l}asz \cite{Haj2}. We take a subsequence $(u_{k_j})_j$ such that
 \begin{equation} \label{p,q quasinorm of ukj minus u}
 ||u_{k_j}-u||_{L^{p,q}(X, \mu)} < 2^{-2j}.
 \end{equation}
 Set $g_j=|u_{k_j}-u|,$ and let $\Gamma \subset \Gamma_{\rm{rect}}$ be the family of curves such that $$\limsup_{j \rightarrow \infty} \int_{\gamma} g_j>0.$$ We want to show that ${\mathrm{Mod}}_{p,q}(\Gamma)=0.$ Denote by $\Gamma_j$ the family of curves in $\Gamma_{\rm{rect}}$ for which $\int_{\gamma} g_j > 2^{-j}.$ Then $2^j g_j \in F(\Gamma_j)$ and hence ${\mathrm{Mod}}_{p,q}(\Gamma_j)<2^{-pj}$ as a consequence of (\ref{p,q quasinorm of ukj minus u}). We notice that $$\Gamma \subset \bigcap_{i=1}^{\infty} \bigcup_{j=i}^{\infty} \Gamma_j.$$
 Thus $${\mathrm{Mod}}_{p,q}(\Gamma)^{1/p} \le \sum_{j=i}^{\infty} {\mathrm{Mod}}_{p,q}(\Gamma_j)^{1/p} \le \sum_{j=i}^{\infty} 2^{-j}=2^{1-i}$$
 for every integer $i \ge 1,$ which implies ${\mathrm{Mod}}_{p,q}(\Gamma)=0.$
 \end{pf}

 \subsection{Upper gradient}

 \begin{Definition} Let $u: X \rightarrow [-\infty, \infty]$ be a Borel function. We say that a Borel function $g: X \rightarrow [0, \infty]$ is an \textit{upper gradient} of $u$ if for \textit{every} rectifiable curve $\gamma$ parametrized by arc length parametrization we have
 \begin{equation} \label{defn of upper gradient}
 |u(\gamma(0))-u(\gamma(\ell(\gamma)))| \le \int_{\gamma} g
 \end{equation}
 whenever both $u(\gamma(0))$ and $u(\gamma(\ell(\gamma)))$ are finite and $\int_{\gamma} g=\infty$ otherwise. We say that $g$ is a $p,q$-weak upper gradient of $u$ if (\ref{defn of upper gradient}) holds on $p,q$-almost every curve $\gamma \in \Gamma_{\rm{rect}}.$

 \end{Definition}

 The weak upper gradients were introduced in the case $p=q$ by Heinonen-Koskela in \cite{HeK}. See also Heinonen \cite{Hei} and Shanmugalingam \cite{Sha1} and \cite{Sha2}.

 If $g$ is an upper gradient of $u$ and $\widetilde{g}=g,$ $\mu$-almost everywhere, is another nonnegative Borel function, then it might happen that $\widetilde{g}$ is not an upper gradient of $u.$ However, we have the following result.

 \begin{Lemma} If $g$ is a $p,q$-weak upper gradient of $u$ and $\widetilde{g}$ is another nonnegative Borel function such that $\widetilde{g}=g$ $\mu$-almost everywhere, then $\widetilde{g}$ is a $p,q$-weak upper gradient of $u.$

 \end{Lemma}

 \begin{pf} Let $\Gamma_1 \subset \Gamma_{\rm{rect}}$ be the family of all nonconstant rectifiable curves $\gamma:[0, \ell(\gamma)] \rightarrow X$ for which $\int_{\gamma} |g-\widetilde{g}|>0.$ The constant sequence $g_n=|g-\widetilde{g}|$ converges to $0$ in $L^{p,q}(X,\mu),$ so from Theorem \ref{uk conv to u in Lpq implies uk conv to u on curves} it follows that ${\mathrm{Mod}}_{p,q}(\Gamma_1)=0$ and $\int_{\gamma} |g-\widetilde{g}|=0$ for every nonconstant rectifiable curve $\gamma: [0, \ell(\gamma)] \rightarrow X$ that is not in $\Gamma_1.$

 Let $\Gamma_2 \subset \Gamma_{\rm{rect}}$ be the family of all nonconstant rectifiable curves $\gamma: [0, \ell(\gamma)] \rightarrow X$ for which the inequality
 $$|u(\gamma(0))-u(\gamma(\ell(\gamma)))| \le \int_{\gamma} g $$ is \textit{not} satisfied. Then ${\mathrm{Mod}}_{p,q}(\Gamma_2)=0.$ Thus ${\mathrm{Mod}}_{p,q}(\Gamma_1 \cup \Gamma_2)=0.$ For every $\gamma \in \Gamma_{\rm{rect}}$ not in $\Gamma_{1} \cup \Gamma_{2}$ we
 have $$ |u(\gamma(0))-u(\gamma(\ell(\gamma)))| \le \int_{\gamma} g= \int_{\gamma} \widetilde{g}.$$
 This finishes the proof.

 \end{pf}

 The next result shows that $p,q$-weak upper gradients can be nicely approximated by upper gradients. The case $p=q$ was proved by Koskela-MacManus \cite{KoM}.

 \begin{Lemma} \label{pq weak upper gradient approximated by upper gradients} If $g$ is a $p,q$-weak upper gradient of $u$ which is finite
 $\mu$-almost everywhere, then for every $\varepsilon>0$ there exists an upper gradient
 $g_{\varepsilon}$ of $u$ such that
 $$ g_{\varepsilon} \ge g \mbox{ everywhere on $X$ and } ||g_{\varepsilon}-g||_{L^{p,q}(X,\mu)} \le \varepsilon.$$

 \end{Lemma}

 \begin{pf} Let $\Gamma \subset \Gamma_{\rm{rect}}$ be the family of all nonconstant
 rectifiable curves $\gamma: [0, \ell(\gamma)] \rightarrow X$ for which the inequality $$|u(\gamma(0))-u(\gamma(\ell(\gamma)))| \le \int_{\gamma} g $$ is \textit{not} satisfied. Then ${\mathrm{Mod}}_{p,q}(\Gamma)=0$ and hence, from Theorem \ref{char of families of curves of pq modulus zero} it follows that there exists $0 \le \rho \in L^{p,q}(X, \mu)$ such that $\int_{\gamma} \rho=\infty$ for every $\gamma \in \Gamma.$ Take $g_{\varepsilon}=g+\varepsilon \rho/||\rho||_{L^{p,q}(X,\mu)}.$ Then $g_{\varepsilon}$ is a nonnegative Borel function and $$|u(\gamma(0))-u(\gamma(\ell(\gamma)))| \le \int_{\gamma} g_{\varepsilon} $$ for \textit{every} curve $\gamma \in \Gamma_{\rm{rect}}.$ This finishes the proof.

 \end{pf}

 If $A$ is a subset of $X$ let $\Gamma_{A}$ be the family of all curves in $\Gamma_{\rm{rect}}$ that intersect $A$ and let $\Gamma_{A}^{+}$ be the family of all curves in $\Gamma_{\rm{rect}}$ such that the Hausdorff one-dimensional measure
 ${\mathcal{H}}_{1}(|\gamma| \cap A)$ is positive. Here and throughout the paper $|\gamma|$ is the image of the curve $\gamma.$

 The following lemma will be useful later in this paper.

 \begin{Lemma} \label{grad for all ui and ui conv to u impl grad for u} Let $u_i: X \rightarrow \Bbb{R},$ $i \ge 1,$ be a sequence of Borel functions such that $g \in L^{p,q}(X)$ is a $p,q$-weak upper gradient for every $u_i, i \ge 1.$ We define $u(x)=\lim_{i \rightarrow \infty} u_i(x)$ and $E=\{ x \in X: |u(x)|=\infty \}.$ Suppose that $\mu(E)=0$ and that $u$ is well-defined outside $E.$ Then $g$ is a $p,q$-weak upper gradient for $u.$

 \end{Lemma}

 \begin{pf} For every $i \ge 1$ we define $\Gamma_{1,i}$ to be the set of all curves $\gamma \in \Gamma_{\rm{rect}}$ for which $$|u_i(\gamma(0))-u_i(\gamma(\ell(\gamma)))| \le \int_{\gamma} g $$ is \textit{not} satisfied. Then ${\mathrm{Mod}}_{p,q}(\Gamma_{1,i})=0$ and hence ${\mathrm{Mod}}_{p,q}(\Gamma_{1})=0,$
 where $\Gamma_{1}=\cup_{i=1}^{\infty} \Gamma_{1,i}.$

 \par Let $\Gamma_0$ be the collection of all paths $\gamma \in \Gamma_{\rm{rect}}$ such that $\int_{\gamma} g=\infty.$ Then we have via Theorem \ref{char of families of curves of pq modulus zero} that ${\mathrm{Mod}}_{p,q}(\Gamma_0)=0$ since $g \in L^{p,q}(X,\mu).$

 Since $\mu(E)=0,$ it follows that ${\mathrm{Mod}}_{p,q}(\Gamma_{E}^{+})=0.$ Indeed, $\infty \cdot \chi_{E} \in F(\Gamma_{E}^{+})$ and $||\infty \cdot \chi_{E}||_{L^{p,q}(X, \mu)}=0.$ Therefore ${\mathrm{Mod}}_{p,q}(\Gamma_0 \cup \Gamma_{E}^{+} \cup \Gamma_1)=0.$

 \par For any path $\gamma$ in the family ${\Gamma_{\rm{rect}}} \setminus (\Gamma_0 \cup \Gamma_{E}^{+} \cup \Gamma_1),$ by the fact that the path
 is not in $\Gamma_{E}^{+}$, there exists a point $y$ in $|\gamma|$ such that
 $y$ is not in $E,$ that is $y \in |\gamma|$ and $|u(y)|<\infty.$ For any point
 $x \in |\gamma|,$ we have (since $\gamma$ is not in $\Gamma_{1,i}$)
 $$|u_{i}(x)|-|u_{i}(y)| \le |u_{i}(x)-u_{i}(y)| \le \int_{\gamma} g<\infty.$$

 Therefore $$|u_{i}(x)| \le |u_{i}(y)| + \int_{\gamma} g.$$ Taking limits on both sides and using the facts that $|u(y)|<\infty$ and that $\gamma$ is not in $\Gamma_0 \cup \Gamma_1,$ we see that
 $$\lim_{i \rightarrow \infty} |u_{i}(x)| \le \lim_{i \rightarrow \infty} |u_{i}(y)| + \int_{\gamma} g=|u(y)|+\int_{\gamma} g<\infty$$ and therefore
 $x$ is not in $E.$ Thus $\Gamma_{E} \subset \Gamma_0 \cup \Gamma_{E}^{+} \cup \Gamma_1$ and ${\mathrm{Mod}}_{p,q}(\Gamma_{E})=0.$

 Next, let $\gamma$ be a path in ${\Gamma_{\rm{rect}}} \setminus (\Gamma_0 \cup \Gamma_{E}^{+} \cup \Gamma_1).$ The above argument showed that
 $|\gamma|$ does not intersect $E.$ If we denote by $x$ and $y$ the endpoints of $\gamma,$ we have
 $$|u(x)-u(y)|=|\lim_{i \rightarrow \infty} u_{i}(x)-\lim_{i \rightarrow \infty} u_{i}(y)|= \lim_{i \rightarrow \infty}  |u_{i}(x)-u_{i}(y)| \le \int_{\gamma} g.$$
 Therefore $g$ is a $p,q$-weak upper gradient for $u$ as well.

 \end{pf}

 The following proposition shows how the upper gradients behave under a change of variable.

 \begin{Proposition} \label{Prop upper grad for F comp u}
  Let $F: \Bbb{R} \rightarrow \Bbb{R}$ be $C^{1}$ and let $u: X \rightarrow \Bbb{R}$ be a Borel function. If $g \in L^{p,q}(X, \mu)$ is a $p,q$-weak upper gradient for $u,$ then $|F'(u)| g$ is a $p,q$-weak upper gradient for $F \circ u.$

 \end{Proposition}

 \begin{pf} Let $\Gamma_{0}$ to be the set of all curves $\gamma \in \Gamma_{\rm{rect}}$ for which $$|u(\gamma(0))-u(\gamma(\ell(\gamma)))| \le \int_{\gamma} g $$ is \textit{not} satisfied. Then ${\mathrm{Mod}}_{p,q}(\Gamma_{0})=0.$ Let $\Gamma_{1} \subset \Gamma_{\rm{rect}}$ be the collection of all curves having a subcurve in $\Gamma_{0}.$ Then $F(\Gamma_{0}) \subset F(\Gamma_{1})$ and hence ${\mathrm{Mod}}_{p,q}(\Gamma_{1}) \le {\mathrm{Mod}}_{p,q}(\Gamma_{0})=0.$

 Let $\Gamma_{2}$ be the set of curves $\gamma \in \Gamma_{\rm{rect}}$ for which
 $\int_{\gamma} g=\infty.$ Then we have via Theorem \ref{char of families of curves of pq modulus zero} that ${\mathrm{Mod}}_{p,q}(\Gamma_2)=0$ since
 $g \in L^{p,q}(X, \mu).$ Thus ${\mathrm{Mod}}_{p,q}(\Gamma_1 \cup \Gamma_2)=0.$

 The claim will follow immediately after we show that
 \begin{equation} \label{upper grad for F comp u and delta}
 |(F \circ u)(\gamma(0))-(F \circ u)(\gamma(\ell(\gamma)))| \le \int_{0}^{\ell(\gamma)} (|F'(u(\gamma(s)))|+\varepsilon) g(\gamma(s)) \, ds.
 \end{equation}
 for all curves $\gamma \in {\Gamma_{\rm{rect}}} \setminus (\Gamma_1 \cup \Gamma_2)$ and for every $\varepsilon>0.$

 So fix $\varepsilon>0$ and choose a curve $\gamma \in {\Gamma_{\rm{rect}}} \setminus (\Gamma_1 \cup \Gamma_2).$ Let $\ell=\ell(\gamma).$ We notice immediately that $u \circ \gamma$ is uniformly continuous on $[0, \ell]$ and $F'$ is uniformly continuous on the compact interval $I:=(u \circ \gamma)([0, \ell]).$ Let $\delta, \delta_{1}>0$ be chosen such that $$|(F'\circ u \circ \gamma)(t)-(F'\circ u \circ \gamma)(s)|+|(u \circ \gamma)(t)-(u \circ \gamma)(s)|<\delta_{1}$$ for all $t,s \in [0, \ell]$ with $|t-s|<\delta$ and such that $$|F'(u)-F'(v)|<\varepsilon \mbox{ for all $u,v \in I$ with $|u-v|<\delta_{1}$}.$$
 Fix an integer $n>1/\delta$ and put $\ell_i=(i\ell)/n, i=0, \ldots, n-1.$
 For every $i=0, \ldots, n-1$ we have
 \begin{eqnarray*}|(F\circ u \circ \gamma)(\ell_{i+1})-(F\circ u \circ \gamma)(\ell_{i})|&=&|F'(t_{i, i+1}))| \, |(u \circ \gamma)(\ell_{i+1})-(u \circ \gamma)(\ell_{i})|\\ &\le& |F'(t_{i, i+1}))| \, \int_{\ell_i}^{\ell_{i+1}} g(\gamma(s)) \, ds,
 \end{eqnarray*}
 where $t_{i, i+1} \in I_{i,i+1}:=(u \circ \gamma)((\ell_{i}, \ell_{i+1})).$
 From the choice of $\delta$ it follows that
 \begin{equation*}|(F\circ u \circ \gamma)(\ell_{i+1})-(F\circ u \circ \gamma)(\ell_{i})| \le \int_{\ell_i}^{\ell_{i+1}} (|F'(u(\gamma(s)))|+\varepsilon) \, g(\gamma(s)) \, ds,
 \end{equation*}
 for every $i=0, \ldots, n-1.$
 If we sum over $i$ we obtain easily (\ref{upper grad for F comp u and delta}). This finishes the proof.
 \end{pf}

 As a direct consequence of Proposition \ref{Prop upper grad for F comp u}, we have the following corollaries.

 \begin{Corollary} \label{gradient of u to power r greater than one} Let $r \in (1, \infty)$ be fixed. Suppose $u: X \rightarrow \Bbb{R}$ is a bounded nonnegative Borel function. If $g \in L^{p,q}(X, \mu)$ is a $p,q$-weak upper gradient of $u,$ then $r u^{r-1} g$ is a $p,q$-weak upper gradient for $u^r.$

 \end{Corollary}

 \begin{pf} Let $M>0$ be such that $0 \le u(x) <M$ for all $x \in X.$ We apply Proposition \ref{Prop upper grad for F comp u} to any $C^1$ function $F: \Bbb{R} \rightarrow \Bbb{R}$ satisfying $F(t)=t^r, 0 \le t \le M.$

 \end{pf}

 \begin{Corollary} \label{gradient of u to power r smaller than one} Let $r \in (0,1)$ be fixed. Suppose that $u: X \rightarrow \Bbb{R}$ is a nonnegative function that has a $p,q$-weak upper gradient $g \in L^{p,q}(X, \mu).$ Then $r (u+\varepsilon)^{r-1} g$ is a $p,q$-weak upper gradient for $(u+\varepsilon)^r$ for all $\varepsilon>0.$

 \end{Corollary}

 \begin{pf} Fix $\varepsilon>0.$ We apply Proposition \ref{Prop upper grad for F comp u} to any $C^1$ function $F: \Bbb{R} \rightarrow \Bbb{R}$ satisfying $F(t)=t^r, \varepsilon \le t <\infty.$

 \end{pf}

 \begin{Corollary} \label{grad for root q of ulq plus u2q} Suppose $1 \le q \le p<\infty.$ Let $u_1, u_2$ be two nonnegative bounded real-valued functions defined on $X.$ Suppose $g_i \in L^{p,q}(X, \mu), i=1,2$ are $p,q$-weak upper gradients for $u_i, i=1,2.$ Then $L^{p,q}(X, \mu) \ni g:=(g_1^q+g_2^q)^{1/q}$ is a $p,q$-weak upper gradient for $u:=(u_1^q+u_2^q)^{1/q}.$

 \end{Corollary}

 \begin{pf} The claim is obvious when $q=1,$ so we assume without loss of generality that $1<q \le p.$ We prove first that $g \in L^{p,q}(X,\mu).$ Indeed, via Remark \ref{identities with powers of norms} it is enough to show that $g^q \in L^{\frac{p}{q}, 1}(X, \mu).$ But $g^q=g_1^q+g_2^q$ and $g_i^q \in L^{\frac{p}{q},1}(X,\mu)$ since $g_i \in L^{p, q}(X, \mu).$ (See Remark \ref{identities with powers of norms}.) This, the fact that $||\cdot||_{L^{\frac{p}{q},1}(X,\mu)}$ is a norm whenever $1<q\le p,$ and another appeal to Remark \ref{identities with powers of norms} yield $g \in L^{p,q}(X,\mu)$ with
 \begin{eqnarray*}
 ||g||_{L^{p,q}(X, \mu)}^{q}&=&||g^q||_{L^{\frac{p}{q},1}(X, \mu)} \le ||g_1^q||_{L^{\frac{p}{q},1}(X, \mu)}+||g_2^q||_{L^{\frac{p}{q},1}(X, \mu)}\\
 &=&||g_1||_{L^{p,q}(X, \mu)}^{q}+||g_2||_{L^{p,q}(X, \mu)}^{q} .
 \end{eqnarray*}

 For $i=1,2$ let $\Gamma_{i,1}$ be the family of nonconstant rectifiable curves in $\Gamma_{\rm{rect}}$ for which
 $$|u_i(x)-u_i(x)| \le \int_{\gamma} g_i$$
 is \textit{not} satisfied. Then ${\mathrm{Mod}}_{p,q}(\Gamma_{i,1})=0$ since $g_i$ is a $p,q$-weak upper gradient for $u_i, i=1,2.$

 Let $\Gamma_{i,2}$ be the family of nonconstant rectifiable curves in $\Gamma_{\rm{rect}}$ for which $\int_{\gamma} g_i=\infty.$ Then for $i=1,2$ we have ${\mathrm{Mod}}_{p,q}(\Gamma_{i,2})=0$ via Theorem \ref{char of families of curves of pq modulus zero} because by hypothesis $g_i \in L^{p,q}(X, \mu), i=1,2.$
 Let $\Gamma_{0}=\Gamma_{1,1} \cup \Gamma_{1,2} \cup \Gamma_{2,1} \cup \Gamma_{2,2}.$ Then ${\mathrm{Mod}}_{p,q}(\Gamma_{0})=0.$

 Fix $\varepsilon>0.$ By applying Corollary \ref{gradient of u to power r greater than one} with $r=q$ to the functions $u_i$ for $i=1,2,$ we see that
 $L^{p,q}(X, \mu) \ni q (u_i+\varepsilon)^{q-1} g_i$ is a $p,q$-weak upper gradient of $(u_i+\varepsilon)^{q}$ for $i=1,2.$ Thus via H\"{o}lder's inequality it follows that $G_{\varepsilon}$ is a $p,q$-weak upper gradient for $U_{\varepsilon},$ where $$G_{\varepsilon}:=q((u_1+\varepsilon)^{q}+(u_2+\varepsilon)^{q})^{(q-1)/q} (g_{1}^q+g_{2}^q)^{1/q} \mbox{ and } U_{\varepsilon}:=(u_1+\varepsilon)^{q}+(u_2+\varepsilon)^{q}.$$
 We notice that $G_{\varepsilon} \in L^{p,q}(X,\mu).$ Indeed,
 $G_{\varepsilon}=q U_{\varepsilon}^{(q-1)/q} g,$ with $U_{\varepsilon}$ nonnegative a bounded and $g \in L^{p,q}(X,\mu),$ so $G_{\varepsilon} \in L^{p,q}(X,\mu).$

 Now we apply Corollary \ref{gradient of u to power r smaller than one} with $r=1/q,$ $u=U_{\varepsilon}$ and $g=G_{\varepsilon}$ to obtain that
 $u_{\varepsilon}:=U_{\varepsilon}^{1/q}$ has $1/q U_{\varepsilon}^{(1-q)/q} G_{\varepsilon}=g$ as a $p,q$-weak upper gradient that belongs to $L^{p,q}(X, \mu).$
 In fact, by looking at the proof of Proposition \ref{Prop upper grad for F comp u},  we see that $$|u_{\varepsilon}(x)-u_{\varepsilon}(y)| \le \int_{\gamma} g$$
 for every curve $\gamma \in \Gamma_{\rm{rect}}$ that is not in $\Gamma_{0}.$
 Letting $\varepsilon \rightarrow 0,$ we obtain the desired conclusion. This finishes the proof of the corollary.

 \end{pf}

 \begin{Lemma} \label{gradient for max of two functions} If $u_i, i=1,2$ are nonnegative real-valued Borel functions in $L^{p,q}(X,\mu)$ with corresponding $p,q$-weak upper gradients $g_i \in L^{p,q}(X,\mu),$ then $g:=\max(g_1, g_2) \in L^{p,q}(X, \mu)$ is a $p,q$-weak upper gradient for $u:=\max(u_1, u_2) \in L^{p,q}(X, \mu).$

 \end{Lemma}

 \begin{pf} It is easy to see that $u, g \in L^{p,q}(X, \mu).$
 For $i=1,2$ let $\Gamma_{0,i} \subset \Gamma_{\rm{rect}}$ be the family of nonconstant rectifiable curves $\gamma$ for which $\int_{\gamma} g_i=\infty.$
 Then we have via Theorem \ref{char of families of curves of pq modulus zero} that ${\mathrm{Mod}}_{p,q}(\Gamma_{0,i})=0$ because $g_i \in L^{p,q}(X, \mu).$ Thus ${\mathrm{Mod}}_{p,q}(\Gamma_{0})=0,$ where
 $\Gamma_{0}=\Gamma_{0,1} \cup \Gamma_{0,2}.$

 For $i=1,2$ let $\Gamma_{i,1} \subset \Gamma_{\rm{rect}}$ be the family of curves $\gamma \in \Gamma_{\rm{rect}} \setminus \Gamma_{0}$ for which $$|u_i(\gamma(0))-u_i(\gamma(\ell(\gamma)))| \le \int_{\gamma} g_i$$ is \textit{not} satisfied. Then ${\mathrm{Mod}}_{p,q}(\Gamma_{1,i})=0$ since $g_i$ is a $p,q$-weak upper gradient for $u_i,$ $i=1,2.$ Thus ${\mathrm{Mod}}_{p,q}(\Gamma_{1})=0,$ where
 $\Gamma_{1}=\Gamma_{1,1} \cup \Gamma_{1,2}.$

 It is easy to see that
 \begin{equation} \label{abs value diff sup le sup abs value diff}
 |u(x)-u(y)| \le \max(|u_1(x)-u_1(y)|, |u_2(x)-u_2(y)|).
 \end{equation}
 On every curve $\gamma \in \Gamma_{\rm{rect}} \setminus (\Gamma_0 \cup \Gamma_1)$
 we have $$|u_i(x)-u_i(y)| \le \int_{\gamma} g_i \le \int_{\gamma} g.$$ This and (\ref{abs value diff sup le sup abs value diff}) show that
 $$|u(x)-u(y)| \le \int_{\gamma} g$$ on every curve $\gamma \in \Gamma_{\rm{rect}} \setminus (\Gamma_0 \cup \Gamma_1).$ This finishes the proof.

 \end{pf}

 \begin{Lemma} \label{gradient truncation above} Suppose $g \in L^{p,q}(X, \mu)$ is a $p,q$-weak upper gradient for $0 \le u \in L^{p,q}(X, \mu).$ Let $\lambda \ge 0$ be fixed. Then  $u_{\lambda}:=\min(u, \lambda) \in L^{p,q}(X, \mu)$ and $g$ is a $p,q$-weak upper gradient for $u_{\lambda}.$

 \end{Lemma}

 \begin{pf} Obviously $0 \le u_{\lambda} \le u$ on $X,$ so it follows via Bennett-Sharpley \cite[Proposition I.1.7]{BS} and Kauhanen-Koskela-Mal{\'{y}} \cite[Proposition 2.1]{KKM} that $u_{\lambda} \in L^{p,q}(X, \mu)$ with $||u_{\lambda}||_{L^{p,q}(X, \mu)} \le ||u||_{L^{p,q}(X, \mu)}.$ The second claim follows immediately since $|u_{\lambda}(x)-u_{\lambda}(y)| \le |u(x)-u(y)|$ for every $x, y \in X.$

 \end{pf}

 \section{Newtonian $L^{p,q}$ spaces}

 We denote by $\widetilde{N}^{1, L^{p,q}}(X,\mu)$ the space of all Borel functions $u \in L^{p,q}(X, \mu)$ that have a $p,q$-weak upper gradient $g \in L^{p,q}(X, \mu).$
 We note that $\widetilde{N}^{1, L^{p,q}}(X,\mu)$ is a vector space, since if $\alpha, \beta \in \Bbb{R}$ and $u_1, u_2 \in \widetilde{N}^{1, L^{p,q}}(X,\mu)$ with respective $p,q$-weak upper gradients $g_1, g_2 \in L^{p,q}(X,\mu),$ then $|\alpha| g_1+ |\beta| g_2 $ is a $p,q$-weak upper gradient of $\alpha u_1+ \beta u_2.$

 \begin{Definition} \label{definition qusinorms N1pq}
 If $u$ is a function in $\widetilde{N}^{1, L^{p,q}}(X,\mu),$ let
 $$||u||_{\widetilde{N}^{1, L^{p,q}}}:=\left\{ \begin{array}{lc}
\left(||u||_{L^{p,q}(X, \mu)}^{q} + \inf_{g} ||g||_{L^{p,q}(X, \mu)}^{q}\right)^{1/q}, & 1 \le q \le p, \\

\left(||u||_{L^{p,q}(X, \mu)}^{p} + \inf_{g} ||g||_{L^{p,q}(X, \mu)}^{p}\right)^{1/p}, &
p<q \le \infty,
\end{array}
\right.
$$
where the infimum is taken over all $p,q$-integrable $p,q$-weak upper gradients of $u.$

 Similarly, let
 $$||u||_{\widetilde{N}^{1, L^{(p,q)}}}:=\left\{ \begin{array}{lc}
\left(||u||_{L^{(p,q)}(X, \mu)}^{q} + \inf_{g} ||g||_{L^{(p,q)}(X, \mu)}^{q}\right)^{1/q}, & 1 \le q \le p, \\

\left(||u||_{L^{(p,q)}(X, \mu)}^{p} + \inf_{g} ||g||_{L^{(p,q)}(X, \mu)}^{p}\right)^{1/p}, &
p<q \le \infty,
\end{array}
\right.
$$
where the infimum is taken over all $p,q$-integrable $p,q$-weak upper gradients of $u.$

 \end{Definition}

 If $u, v$ are functions in $\widetilde{N}^{1, L^{p,q}}(X, \mu),$ let $u \sim v$ if $||u-v||_{\widetilde{N}^{1, L^{p,q}}}=0.$ It is easy to see that $\sim$ is an equivalence relation that partitions $\widetilde{N}^{1, L^{p,q}}(X, \mu)$ into equivalence classes. We define the space $N^{1, L^{p,q}}(X, \mu)$ as the quotient
 $\widetilde{N}^{1, L^{p,q}}(X, \mu)/{\sim}$ and
 $$||u||_{N^{1, L^{p,q}}}=||u||_{\widetilde{N}^{1, L^{p,q}}} \mbox{ and } ||u||_{N^{1, L^{(p,q)}}}=||u||_{\widetilde{N}^{1, L^{(p,q)}}}$$

 \begin{Remark} \label{comparison of Newtonian norms and quasinorms} Via Lemma \ref{pq weak upper gradient approximated by upper gradients} and Corollary \ref{Convergence of the quasinorm of the sequence}, it is easy to see that the infima in Definition \ref{definition qusinorms N1pq} could as well be taken over all $p,q$-integrable upper gradients of $u.$ We also notice (see the discussion before Definition \ref{definition of abs continuous norm}) that $||\cdot||_{N^{1, L^{(p,q)}}}$ is a norm whenever $1<p<\infty$ and $1 \le q \le \infty$ , while
 $||\cdot||_{N^{1, L^{p,q}}}$ is a norm when $1 \le q \le p<\infty.$ Moreover
 (see the discussion before Definition \ref{definition of abs continuous norm}) $$||u||_{N^{1, L^{p,q}}} \le ||u||_{N^{1, L^{(p,q)}}} \le p' ||u||_{N^{1, L^{p,q}}}$$
 for every $1<p<\infty,$ $1 \le q \le \infty$ and $u \in N^{1, L^{p,q}}(X, \mu).$
 \end{Remark}

 \begin{Definition} Let $u: X \rightarrow [-\infty, \infty]$ be a given function. We say that
 \par (i) $u$ is absolutely continuous along a rectifiable curve $\gamma$ if $u \circ \gamma$ is absolutely continuous on $[0, \ell(\gamma)].$

 \par (ii) $u$ is absolutely continuous on $p,q$-almost every curve (has $ACC_{p,q}$ property) if for $p,q$-almost every $\gamma \in \Gamma_{\rm{rect}},$ $u \circ \gamma$ is absolutely continuous.
 \end{Definition}

 \begin{Proposition} \label{N1pq is a subset of ACCpq} If $u$ is a function in $\widetilde{N}^{1, L^{p,q}}(X,\mu),$ then $u$ is $ACC_{p,q}.$

 \end{Proposition}

 \begin{pf} We follow Shanmugalingam \cite{Sha1}. By the definition of $\widetilde{N}^{1, L^{p,q}}(X,\mu),$ $u$ has a $p,q$-weak upper gradient $g \in L^{p,q}(X, \mu).$ Let $\Gamma_{0}$ be the collection of all curves in $\Gamma_{\rm{rect}}$ for which $$|u(\gamma(0))-u(\gamma(\ell(\gamma)))| \le \int_{\gamma} g$$ is \textit{not} satisfied.
 Then by the definition of $p,q$-weak upper gradients, ${\mathrm{Mod}}_{p,q}(\Gamma_{0})=0.$ Let $\Gamma_{1}$ be the collection of all curves in $\Gamma_{\rm{rect}}$ that have a subcurve in $\Gamma_{0}.$ Then ${\mathrm{Mod}}_{p,q}(\Gamma_{1}) \le {\mathrm{Mod}}_{p,q}(\Gamma_{0})=0.$

 Let $\Gamma_{2}$ be the collection of all curves in $\Gamma_{\rm{rect}}$ such that
 $\int_{\gamma} g=\infty.$ Then ${\mathrm{Mod}}_{p,q}(\Gamma_{2})=0$ because $g \in L^{p,q}(X, \mu).$ Hence ${\mathrm{Mod}}_{p,q}(\Gamma_{1} \cup \Gamma_{2})=0.$ If $\gamma$ is a curve in $\Gamma_{\rm{rect}} \setminus (\Gamma_{1} \cup \Gamma_{2}),$
 then $\gamma$ has no subcurves in $\Gamma_{0},$ and hence
 $$|u(\gamma(\beta))-u(\gamma(\alpha))| \le \int_{\alpha}^{\beta} g(\gamma(t)) \, dt, \mbox{ provided } [\alpha, \beta] \subset [0, \ell(\gamma)].$$
 This implies the absolute continuity of $u \circ \gamma$ as a consequence of the absolute continuity of the integral. Therefore $u$ is absolutely continuous on every curve $\gamma$ in $\Gamma_{\rm{rect}} \setminus (\Gamma_{1} \cup \Gamma_{2}).$

 \end{pf}

 \begin{Lemma} Suppose $u$ is a function in $\widetilde{N}^{1, L^{p,q}}(X, \mu)$ such that $||u||_{L^{p,q}(X, \mu)}=0.$ Then the family
 $$\Gamma=\{ \gamma \in \Gamma_{\rm{rect}}: u(x) \neq 0 \mbox{ for some } x \in |\gamma|  \} $$
 has zero $p,q$-modulus.

 \end{Lemma}

 \begin{pf} We follow Shanmugalingam \cite{Sha1}. Since $||u||_{L^{p,q}(X,\mu)}=0,$ the set $E=\{ x \in X: u(x) \neq 0 \}$ has measure zero.
 With the notation introduced earlier, we have
 $$\Gamma=\Gamma_{E}=\Gamma_{E}^{+} \cup (\Gamma_{E} \setminus \Gamma_{E}^{+}).$$
 We can disregard the family $\Gamma_{E}^{+},$ since
 $${\mathrm{Mod}}_{p,q}(\Gamma_{E}^{+}) \le ||\infty \cdot \chi_{E}||_{L^{p,q}(X, \mu)}^{p}=0,$$
 where $\chi_{E}$ is the characteristic function of the set $E.$ The curves $\gamma$ in $\Gamma_{E} \setminus \Gamma_{E}^{+}$ intersect $E$ only on a set of linear measure zero, and hence with respect to the linear measure almost everywhere on $\gamma$ the function $u$ is equal to zero. Since $\gamma$ also intersects $E,$ it follows that
 $u$ is \textit{not} absolutely continuous on $\gamma.$ By Proposition
 \ref{N1pq is a subset of ACCpq}, we have ${\mathrm{Mod}}_{p,q}(\Gamma_{E} \setminus \Gamma_{E}^{+})=0,$ yielding ${\mathrm{Mod}}_{p,q}(\Gamma)=0.$ This finishes the proof.

 \end{pf}

 \begin{Lemma} \label{gradient trimming on closed sets} Let $F$ be a closed subset of $X.$ Suppose that $u: X \rightarrow [-\infty, \infty]$ is a Borel $ACC_{p,q}$ function that is constant $\mu$-almost everywhere on $F.$ If $g \in L^{p,q}(X,\mu)$
 is a $p,q$-weak upper gradient of $u,$ then $g \chi_{X \setminus F}$ is a $p,q$-weak upper gradient of $u.$

 \end{Lemma}

 \begin{pf} We can assume without loss of generality that $u=0$ $\mu$-almost everywhere on $F.$ Let $E=\{ x \in F: u(x) \neq 0 \}.$ Then by assumption $\mu(E)=0.$ Hence ${\mathrm{Mod}}_{p,q}(\Gamma_{E}^{+})=0$ because $\infty \cdot \chi_{E} \in F(\Gamma_{E}^{+}).$

 Let $\Gamma_{0} \subset \Gamma_{\rm{rect}}$ be the family of curves on which $u$ is not absolutely continuous or on which
 $$|u(\gamma(0))-u(\gamma(\ell(\gamma)))| \le \int_{\gamma} g$$
 is \textit{not} satisfied. Then ${\mathrm{Mod}}_{p,q}(\Gamma_{0})=0.$ Let $\Gamma_{1} \subset \Gamma_{\rm{rect}}$ be the family of curves that have a subcurve in $\Gamma_{0}.$ Then $F(\Gamma_{0}) \subset F(\Gamma_{1})$ and thus ${\mathrm{Mod}}_{p,q}(\Gamma_{1}) \le {\mathrm{Mod}}_{p,q}(\Gamma_{0})=0.$

 Let $\Gamma_{2} \subset \Gamma_{\rm{rect}}$ be the family of curves on which $\int_{\gamma} g=\infty.$ Then via Theorem
 \ref{char of families of curves of pq modulus zero} we have ${\mathrm{Mod}}_{p,q}(\Gamma_{2})=0$ because
 $g \in L^{p,q}(X, \mu).$

 Let $\gamma:[0, \ell(\gamma)] \rightarrow X$ be a curve in $\Gamma_{\rm{rect}} \setminus (\Gamma_{1} \cup \Gamma_{2} \cup \Gamma_{E}^{+})$ connecting $x$ and $y.$ We show that
 $$|u(x)-u(y)| \le \int_{\gamma} g \chi_{X \setminus F} $$ for every such curve $\gamma.$

 The cases $|\gamma| \subset F \setminus E$ and $|\gamma| \subset (X \setminus F) \cup E$ are trivial. So is the case when both $x$ and $y$ are in $F \setminus E.$
 Let $K:= (u \circ \gamma)^{-1}(\{0\}).$ Then $K$ is a compact subset of $[0, \ell(\gamma)]$ because $u \circ \gamma$ is continuous on $[0, \ell(\gamma)].$ Hence $K$ contains its lower bound $c$ and its upper bound $d.$ Let $x_1=\gamma(c)$ and $y_1=\gamma(d).$

 Suppose that both $x$ and $y$ are in $(X \setminus F) \cup E.$ Then we see that $[c,d] \subset (0, \ell(\gamma))$ and $\gamma([0,c) \cup (d, \ell(\gamma)]) \subset (X \setminus F) \cup E.$

 Moreover, since $\gamma$ is not in $\Gamma_{1}$ and $u(x_1)=u(y_1),$ we have
 $$|u(x)-u(y)| \le |u(x)-u(x_1)|+|u(y_1)-u(y)| \le \int_{\gamma([0,c])} g + \int_{\gamma([d, \ell(\gamma)])} g \le \int_{\gamma} g \chi_{X \setminus F}$$
 because the subcurves $\gamma|_{[0,c]}$ and $\gamma|_{[d. \ell(\gamma)]}$
 intersect $E$ on a set of Hausdorff $1$-measure zero.

 Suppose now by symmetry that $x \in (X \setminus F) \cup E$ and $y \in F \setminus E.$ This means in terms of our notation that $c>0$ and $d=\ell(\gamma).$ We notice that $\gamma([0,c)) \subset (X \setminus F) \cup E$ and $u(x_1)=u(y)$ and thus
 $$|u(x)-u(y)|=|u(x)-u(x_1)| \le \int_{\gamma([0,c])} g \le \int_{\gamma} g \chi_{X \setminus F}$$
 because the subcurve $\gamma|_{[0,c]}$ intersect $E$ on a set of Hausdorff $1$-measure zero.

 This finishes the proof of the lemma.

 \end{pf}

 \begin{Lemma} \label{two pq weak upper gradients and a closed set} Assume that $u \in N^{1, L^{p,q}}(X, \mu),$ and that $g,h \in L^{p,q}(X, \mu)$ are $p,q$-weak upper gradients of $u.$ If $F \subset X$ is a closed set, then
 $$\rho=g \chi_{F} + h \chi_{X \setminus F}$$ is a $p,q$-weak upper gradient of $u$ as well.

 \end{Lemma}

 \begin{pf} We follow Haj{\l}asz \cite{Haj2}. Let $\Gamma_{1} \subset \Gamma_{\rm{rect}}$ be the family of curves on which $\int_{\gamma} (g+h) =\infty.$ Then via Theorem
 \ref{char of families of curves of pq modulus zero} it follows that ${\mathrm{Mod}}_{p,q}(\Gamma_{1})=0$ because $g+h \in L^{p,q}(X, \mu).$

 Let $\Gamma_{2} \subset \Gamma_{\rm{rect}}$ be the family of curves on which $u$ is \textit{not} absolutely continuous. Then via Proposition
 \ref{N1pq is a subset of ACCpq} we see that ${\mathrm{Mod}}_{p,q}(\Gamma_2)=0.$

 Let $\Gamma_{3}^{'} \subset \Gamma_{\rm{rect}}$ be the family of curves on which
 $$|u(\gamma(0))-u(\gamma(\ell(\gamma)))| \le \min \left(\int_{\gamma} g, \int_{\gamma} h \right)$$
 is \textit{not} satisfied. Let $\Gamma_{3} \subset \Gamma_{\rm{rect}}$ be the family of curves which contain subcurves belonging to $\Gamma_{3}^{'}.$ Since $F(\Gamma_{3}^{'}) \subset F(\Gamma_{3}),$ we have ${\mathrm{Mod}}_{p,q}(\Gamma_{3}) \le {\mathrm{Mod}}_{p,q}(\Gamma_{3}^{'})=0.$ Now it remains to show that
 $$|u(\gamma(0))-u(\gamma(\ell(\gamma)))| \le \int_{\gamma} \rho $$
 for all $\gamma \in \Gamma_{\rm{rect}} \setminus (\Gamma_{1} \cup \Gamma_{2} \cup \Gamma_{3}).$
 If $|\gamma| \subset F$ or $|\gamma| \subset X \setminus F,$ then the inequality is obvious. Thus we can assume that the image $|\gamma|$ has a nonempty intersection both with $F$ and with $X \setminus F.$

 The set $\gamma^{-1}(X \setminus F)$ is open and hence it consists of a countable (or finite) number of open and disjoint intervals. Assume without loss of generality that there are countably many such intervals. Denote these intervals by $((t_i, s_i))_{i=1}^{\infty}.$ Let $\gamma_{i}=\gamma|_{[t_i, s_i]}.$ We have
 \begin{eqnarray*}
 |u(\gamma(0))-u(\gamma(\ell(\gamma)))| &\le& |u(\gamma(0))-u(\gamma(t_1))|+|u(\gamma(t_1))-u(\gamma(s_1))|\\
 &+&|u(\gamma(s_1))-u(\gamma(\ell(\gamma)))| \le \int_{\gamma \setminus \gamma_{1}} g + \int_{\gamma_1} h,
 \end{eqnarray*}
 where $\gamma \setminus \gamma_{1}$ denotes the two curves obtained from $\gamma$ by removing the interior part $\gamma_{1},$ that is the curves $\gamma|_{[0, t_1]}$ and $\gamma|_{[s_1,b]}.$ Similarly we can remove a larger number of subcurves of $\gamma.$
 This yields
 $$|u(\gamma(0))-u(\gamma(\ell(\gamma)))| \le \int_{\gamma \setminus \cup_{i=1}^{n} \gamma_i} g + \int_{\cup_{i=1}^{n} \gamma_i} h$$
 for each positive integer $n.$ By applying Lebesgue dominated convergence theorem to the curve integral on $\gamma,$ we obtain
 $$|u(\gamma(0))-u(\gamma(\ell(\gamma)))| \le \int_{\gamma} g \chi_{F} + \int_{\gamma} h \chi_{X \setminus F}= \int_{\gamma} \rho. $$

 \end{pf}

 Next we show that when $1<p<\infty$ and $1 \le q<\infty,$ every function $u \in N^{1,L^{p,q}}(X, \mu)$ has a `smallest' $p,q$-weak upper gradient. For the case $p=q$ see Kallunki-Shanmugalingam \cite{KaS} and Shanmugalingam \cite{Sha2}.

 \begin{Theorem} \label{minimal pq weak upper gradient when q is finite}
 Suppose that $1<p<\infty$ and $1 \le q<\infty.$ For every $u \in N^{1,L^{p,q}}(X, \mu),$ there exists the least $p,q$-weak upper gradient $g_{u} \in L^{p,q}(X,\mu)$ of $u.$ It is smallest in the sense that if $g \in L^{p,q}(X,\mu)$ is another $p,q$-weak upper gradient of $u,$ then $g \ge g_{u}$ $\mu$-almost everywhere.

 \end{Theorem}

 \begin{pf} We follow Haj{\l}asz \cite{Haj2}. Let $m=\inf_{g} ||g||_{L^{p,q}(X,\mu)},$ where the infimum is taken over the set of all $p,q$-weak upper gradients of $u.$ It suffices to show that there exists a $p,q$-weak upper gradient $g_{u}$ of $u$ such that $||g_{u}||_{L^{p,q}(X, \mu)}=m.$ Indeed, if we suppose that $g \in L^{p,q}(X,\mu)$ is another $p,q$-weak upper gradient of $u$ such that the set $\{ g<g_{u} \}$ has positive measure, then by the inner regularity of the measure $\mu$ there exists a closed set $F \subset \{ g<g_{u} \}$ such that $\mu(F)>0.$ Via Lemma
 \ref{two pq weak upper gradients and a closed set} it follows that the function $\rho:=g \chi_{F} + g_{u} \chi_{X \setminus F}$ is a $p,q$-weak upper gradient. Via Kauhanen-Koskela-Mal{\'{y}} \cite[Proposition 2.1]{KKM} that would give $||\rho||_{L^{p,q}(X,\mu)}<||g_{u}||_{L^{p,q}(X,\mu)}=m,$ in contradiction with the minimality of $||g_{u}||_{L^{p,q}(X,\mu)}.$

 Thus it remains to prove the existence of a $p,q$-weak upper gradient $g_{u}$ such that $||g_{u}||_{L^{p,q}(X, \mu)}=m.$ Let $(g_i)_{i=1}^{\infty}$ be a sequence of $p,q$-weak upper gradients of $u$ such that $||g_{i}||_{L^{p,q}(X, \mu)}<m+2^{-i}.$ We will show that it is possible to modify the sequence $(g_i)$ in such a way that we will obtain a new sequence of $p,q$-weak upper gradients $(\rho_{i})$ of $u$ satisfying
 $$||\rho_{i}||_{L^{p,q}(X, \mu)}<m+2^{1-i}, \quad \rho_{1} \ge \rho_{2} \ge \rho_{3} \ge \ldots \mbox{ $\mu$-almost everywhere}.$$
 The sequence $(\rho_{i})_{i=1}^{\infty}$ will be defined by induction. We set $\rho_{1}=g_{1}.$ Suppose the $p,q$-weak upper gradients $\rho_{1}, \rho_{2}, \ldots, \rho_{i}$ have already been chosen. We will now define $\rho_{i+1}.$ Since $\rho_{i} \in L^{p,q}(X,\mu),$ the measure $\mu$ is inner regular and the $(p,q)$-norm has the absolute continuity property whenever $1<p<\infty$ and $1 \le q<\infty$ (see the discussion after Definition \ref{definition of abs continuous norm}), there exists a closed set $F \subset \{ g_{i+1}<\rho_{i} \}$ such that
 $$||\rho_{i} \chi_{ \{ g_{i+1}<\rho_{i} \} \setminus F }||_{L^{p,q}(X, \mu)}<2^{-i-1}.$$
 Now we set $\rho_{i+1}=g_{i+1} \chi_{F} + \rho_{i} \chi_{X \setminus F}.$ Then $$\rho_{i+1} \le \rho_{i} \mbox { and }
 \rho_{i+1} \le g_{i+1} \chi_{F \cup \{ g_{i+1} \ge \rho_{i} \} }+ \rho_{i} \chi_{ \{ g_{i+1}<\rho_{i} \} \setminus F }.$$
 Suppose first that $1 \le q \le p.$ Since $||\cdot||_{L^{p,q}(X,\mu)}$ is a norm, we see that
 \begin{eqnarray*}
 ||\rho_{i+1}||_{L^{p,q}(X, \mu)} &\le& ||g_{i+1} \chi_{F \cup \{ g_{i+1} \ge \rho_{i} \} }||_{L^{p,q}(X, \mu)}+ ||\rho_{i} \chi_{ \{ g_{i+1}<\rho_{i} \} \setminus F }||_{L^{p,q}(X, \mu)}\\
 &<& m+2^{-i-1}+2^{-i-1}=m+2^{-i}.
 \end{eqnarray*}
 Suppose now that $p<q<\infty.$ Then we have via Proposition \ref{subadd p le q}
 \begin{eqnarray*}
 ||\rho_{i+1}||_{L^{p,q}(X, \mu)}^{p} &\le& ||g_{i+1} \chi_{F \cup \{ g_{i+1} \ge \rho_{i} \} }||_{L^{p,q}(X, \mu)}^{p}+ ||\rho_{i} \chi_{ \{ g_{i+1}<\rho_{i} \} \setminus F }||_{L^{p,q}(X, \mu)}^{p}\\
 &<& (m+2^{-i-1})^{p}+2^{-p(i+1)} < (m+2^{-i})^{p}.
 \end{eqnarray*}
 Thus, no matter what $q \in [1, \infty)$ is, we showed that $m \le ||\rho_{i+1}||_{L^{p,q}(X, \mu)}<m+2^{-i}.$ The sequence of $p,q$-weak upper gradients
 $(\rho_i)_{i=1}^{\infty}$ converges pointwise to a function $\rho.$ The absolute continuity of the $(p,q)$-norm (see Bennett-Sharpley \cite[Proposition I.3.6]{BS} and the discussion after Definition \ref{definition of abs continuous norm}) yields $$\lim_{i \rightarrow \infty} ||\rho_{i}-\rho||_{L^{p,q}(X, \mu)}=0.$$ Obviously $||\rho||_{L^{p,q}(X,\mu)}=m.$
 The proof will be finished as soon as we show that $\rho$ is a $p,q$-weak upper gradient for $u.$

 By taking a subsequence if necessary, we can assume that $||\rho_{i}-\rho||_{L^{p,q}(X, \mu)} \le 2^{-2i}$ for every $i \ge 1.$

 Let $\Gamma_{1} \subset \Gamma_{\rm{rect}}$ be the family of curves on which $\int_{\gamma} (\rho+\rho_{i}) =\infty$ for some $i \ge 1.$ Then via Theorem \ref{char of families of curves of pq modulus zero} and the subadditivity of ${\mathrm{Mod}}_{p,q}(\cdot)^{1/p}$ we see that ${\mathrm{Mod}}_{p,q}(\Gamma_{1})=0$
 since $\rho+\rho_{i} \in L^{p,q}(X, \mu)$ for every $i \ge 1.$

 For any integer $i \ge 1$ let $\Gamma_{2,i} \subset \Gamma_{\rm{rect}}$ be the family of curves for which $$|u(\gamma(0))-u(\gamma(\ell(\gamma)))| \le \int_{\gamma} \rho_{i}$$ is \textit{not} satisfied. Then ${\mathrm{Mod}}_{p,q}(\Gamma_{2,i})=0$ because $\rho_{i}$ is a $p,q$-weak upper gradient for $u.$ Let $\Gamma_{2}=\cup_{i=1}^{\infty} \Gamma_{2,i}.$

 Let $\Gamma_{3} \subset \Gamma_{\rm{rect}}$ be the family of curves for which
 $\limsup_{i \rightarrow \infty} \int_{\gamma} |\rho_{i}-\rho|>0.$ Then it follows via Theorem \ref{uk conv to u in Lpq implies uk conv to u on curves} that ${\mathrm{Mod}}_{p,q}(\Gamma_{3})=0.$

 Let $\gamma$ be a curve in $\Gamma_{\rm{rect}} \setminus (\Gamma_{1} \cup \Gamma_{2} \cup \Gamma_{3}).$ On any such curve we have (since $\gamma$ is not in $\Gamma_{2,i}$)
 $$|u(\gamma(0))-u(\gamma(\ell(\gamma)))| \le \int_{\gamma} \rho_{i} \mbox{ for every $i \ge 1$}.$$
 By letting $i \rightarrow \infty,$ we obtain (since $\gamma$ is not in $\Gamma_{1} \cup \Gamma_{3}$)
 $$|u(\gamma(0))-u(\gamma(\ell(\gamma)))| \le \lim_{i \rightarrow \infty} \int_{\gamma} \rho_{i} = \int_{\gamma} \rho<\infty.$$
 This finishes the proof of the theorem.

 \end{pf}

 \section{Sobolev $p,q$-capacity}

 In this section, we establish a general theory of the Sobolev-Lorentz $p,q$-capacity in
 metric measure spaces. If $(X,d,\mu)$ is a metric measure space,
 then the Sobolev $p,q$-capacity of a set $E \subset X$ is
 $${\mathrm{Cap}}_{p,q}(E)= \inf \{||u||_{N^{1,L^{p,q}}}^p : u \in {\mathcal{A}}(E) \},$$
 where $${\mathcal{A}}(E)=\{u \in N^{1,L^{p,q}}(X, \mu): u \ge 1  \mbox { on E}  \}.$$

 We call ${\mathcal{A}}(E)$ the set of \textit{admissible functions} for $E.$ If ${\mathcal{A}}(E)=\emptyset,$ then ${\mathrm{Cap}}_{p,q}(E)=\infty.$

 \begin{Remark} It is easy to see that we can consider only admissible functions $u$ for which $0 \le u \le 1.$ Indeed, for $u \in {\mathcal{A}}(E),$ let $v:=\min(u_{+}, 1),$ where $u_{+}=\max(u, 0).$ We notice that $|v(x)-v(y)| \le |u(x)-u(y)|$ for every $x,y$ in $X,$ which implies that every $p,q$-weak upper gradient for $u$ is also a $p,q$-weak upper gradient for $v.$ This implies that $v \in {\mathcal{A}}(E)$ and $||v||_{N^{1,L^{p,q}}} \le ||u||_{N^{1,L^{p,q}}}.$

 \end{Remark}

 \subsection{Basic properties of the Sobolev $p,q$-capacity}

 A capacity is a monotone, subadditive set function. The following theorem expresses,
 among other things, that this is true for the Sobolev $p,q$-capacity.

 \begin{Theorem}\label{p,q Cap Thm} Suppose that $1<p<\infty$ and $1 \le q \le \infty.$ Suppose that $(X,d,\mu)$ is a complete metric measure space. The set function
 $E \mapsto {\mathrm{Cap}}_{p,q}(E),$ $E \subset X,$ enjoys the following properties:

 \par {\rm{(i)}}  If $E_{1} \subset E_{2},$ then ${\mathrm{Cap}}_{p,q}(E_{1})
 \le {\mathrm{Cap}}_{p,q}(E_{2}).$

 \par {\rm{(ii)}} Suppose that $\mu$ is nonatomic. Suppose that $1<q \le p.$ If $E_{1} \subset E_{2} \subset \ldots \subset E=\bigcup_{i=1}^{\infty} E_{i} \subset X,$ then $${\mathrm{Cap}}_{p,q}(E)=\lim_{i \rightarrow \infty} {\mathrm{Cap}}_{p,q}(E_{i}).$$

 \par {\rm{(iii)}} Suppose that $p<q \le \infty.$ If $E=\bigcup_{i=1}^{\infty} E_{i} \subset X,$ then $${\mathrm{Cap}}_{p,q}(E) \le \sum_{i=1}^{\infty} {\mathrm{Cap}}_{p,q}(E_{i}).$$

 \par {\rm{(iv)}} Suppose that $1 \le q \le p.$ If $E=\bigcup_{i=1}^{\infty} E_{i} \subset X,$ then $${\mathrm{Cap}}_{p,q}(E)^{q/p} \le \sum_{i=1}^{\infty} {\mathrm{Cap}}_{p,q}(E_{i})^{q/p}.$$

 \end{Theorem}

 \begin{pf} Property (i) is an immediate consequence of the definition.

 \vskip 2mm

 (ii) Monotonicity yields
 \begin{equation*}
 L:=\lim_{i \rightarrow \infty} {\mathrm{Cap}}_{p,q}(E_{i}) \le {\mathrm{Cap}}_{p,q}(E).
 \end{equation*}
 To prove the opposite inequality, we may assume without loss of generality that $L< \infty.$
 The reflexivity of $L^{p,q}(X,\mu)$ (guaranteed by the nonatomicity of $\mu$ whenever $1<q \le p<\infty$) will be used here in order to prove the opposite inequality.

 Let $\varepsilon >0$ be fixed. For every $i=1,2, \ldots$ we choose $u_{i} \in
 {\mathcal{A}}(E_{i}),$ $0 \le u_{i} \le 1$ and a corresponding upper gradient $g_{i}$
 such that
\begin{equation}\label{Bv 1}
 ||u_i||_{N^{1,L^{p,q}}}^{q}< {\mathrm{Cap}}_{p,q}(E_{i})^{q/p}+ \varepsilon \le L^{q/p}+ \varepsilon.
\end{equation}
 We notice that $u_{i}$ is a bounded sequence in $N^{1,L^{p,q}}(X, \mu).$ Hence there exists a subsequence, which we denote again by $u_{i}$ and functions $u, g \in L^{p,q}(X, \mu)$
 such that $u_{i} \rightarrow u$ weakly in $L^{p,q}(X, \mu)$ and $g_{i} \rightarrow g$ weakly in $L^{p,q}(X, \mu).$ It is easy to see that
 $$u \ge 0 \mbox{ $\mu$-almost everywhere and } g \ge 0 \mbox{ $\mu$-almost everywhere}.$$
 Indeed, since $u_i$ converges weakly to $u$ in $L^{p,q}(X,\mu)$ which is the dual of $L^{p',q'}(X, \mu)$ (see Hunt \cite[p.\ 262]{Hu}), we have
 $$\lim_{i \rightarrow \infty} \int_{X} u_i(x) \varphi(x) \, d \mu(x) = \int_{X} u(x) \varphi(x) \, d\mu(x)$$
 for all $\varphi \in L^{p',q'}(X, \mu).$ For nonnegative functions $\varphi \in L^{p',q'}(X, \mu)$ this yields
 $$0 \le \lim_{i \rightarrow \infty} \int_{X} u_i(x) \varphi(x) \, d \mu(x)=\int_{X} u(x) \varphi(x) \, d\mu(x),$$ which easily implies $u \ge 0$ $\mu$-almost everywhere on $X.$ Similarly we have $g \ge 0$ $\mu$-almost everywhere on $X.$

 From the weak-$*$ lower semicontinuity of the $p,q$-norm (see Bennett-Sharpley \cite[Proposition II.4.2, Definition IV.4.1 and Theorem IV.4.3]{BS} and Hunt
 \cite[p.\ 262]{Hu}), it follows that
 \begin{eqnarray} \label{Fatou for the pq quasinorm}
 & & ||u||_{L^{p,q}(X, \mu)} \le \liminf_{i \rightarrow \infty} ||u_i||_{L^{p,q}(X, \mu)} \mbox{  and  } ||g||_{L^{p,q}(X, \mu)} \le \liminf_{i \rightarrow \infty} ||g_i||_{L^{p,q}(X, \mu)}.
 \end{eqnarray}

 Using Mazur's lemma simultaneously for $u_i$ and $g_{i},$ we obtain sequences $v_i$ with correspondent upper gradients $\widetilde{g}_i$
 such that $v_{i} \in {\mathcal{A}}(E_{i}),$ $v_{i} \rightarrow u$ in $L^{p,q}(X, \mu)$ and $\mu$-almost everywhere and
 $\widetilde{g}_{i} \rightarrow g$ in $L^{p,q}(X,\mu)$ and $\mu$-almost everywhere.
 These sequences can be found in the following way. Let $i_{0}$ be fixed. Since every subsequence of $(u_{i}, g_{i})$
 converges to $(u, g)$ weakly in the reflexive space $L^{p,q}(X, \mu) \times L^{p,q}(X, \mu),$ we may use the Mazur lemma (see Yosida \cite[p.\ 120]{Yos}) for the subsequence $(u_{i}, g_i), i \ge i_{0}.$

 We obtain finite convex combinations $v_{i_0}$ and $\widetilde{g}_{i_0}$ of the functions $u_{i}$ and $g_{i},$ $i\ge i_{0}$
 as close as we want in $L^{p,q}(X, \mu)$ to $u$ and $g$ respectively. For every $i=i_{0}, i_{0}+1, \ldots$
 we see that $u_{i}=1$ in
 $E_{i} \supset E_{i_0}.$ The intersection of finitely many supersets of $E_{i_0}$ contains $E_{{i_0}}.$ Therefore, $v_{i_0}$ equals $1$ on $E_{i_0}.$ It is easy to see that $\widetilde{g}_{i_0}$ is an upper gradient for $v_{i_0}.$
 Passing to subsequences if necessary, we
 may assume that $v_i$ converges to $u$ pointwise $\mu$-almost everywhere, that $\widetilde{g}_i$ converges to $g$ pointwise $\mu$-almost everywhere and that for every $i=1,2, \ldots$ we have
 \begin{equation}
 \label{Bv 3} ||v_{i+1}-v_{i}||_{L^{p,q}(X, \mu)} + ||\widetilde{g}_{{i+1}}-\widetilde{g}_{i}||_{L^{p,q}(X, \mu)} \le 2^{-i}.
 \end{equation}

 Since $v_i$ converges to $u$ in $L^{p,q}(X, \mu)$ and pointwise $\mu$-almost everywhere on $X$ while $\widetilde{g}_i$ converges to $g$ in $L^{p,q}(X, \mu)$ and pointwise $\mu$-almost everywhere on $X$ it follows via Corollary \ref{Convergence of the quasinorm of the sequence} that
 \begin{eqnarray}
 & &\lim_{i \rightarrow \infty} ||v_i||_{L^{p,q}(X, \mu)} =||u||_{L^{p,q}(X, \mu)} \mbox{ and } \lim_{i \rightarrow \infty} ||\widetilde{g}_i||_{L^{p,q}(X, \mu)} =||g||_{L^{p,q}(X, \mu)}.
 \end{eqnarray}
 This, (\ref{Bv 1}) and (\ref{Fatou for the pq quasinorm}) yield
 \begin{equation} \label{norm of (u,g) controlled by norm of L to some power}
 ||u||_{L^{p,q}(X, \mu)}^{q}+||g||_{L^{p,q}(X, \mu)}^{q}= \lim_{i \rightarrow \infty} ||v_i||_{N^{1,L^{p,q}}}^{q} \le L^{q/p}+\varepsilon.
 \end{equation}

 For $j=1,2, \ldots$ we set $$w_{j}=\sup_{i \ge j} v_{i} \mbox{ and } \widehat{g}_j=\sup_{i \ge j} \widetilde{g}_{i}.$$ It is easy to see that $w_j=1$ on $E.$ We claim that $\widehat{g}_j$ is a $p,q$-weak upper gradient for $w_j.$ Indeed, for every $k>j,$ let
 $$w_{j,k}=\sup_{k \ge i \ge j} v_{i}.$$
 Via Lemma \ref{gradient for max of two functions} and finite induction, it follows easily that $\widehat{g}_j$ is a $p,q$-weak upper gradient for every $w_{j,k}$ whenever $k>j.$ It is easy to see that $w_j= \lim_{k\rightarrow \infty} w_{j,k}$ pointwise in $X.$ This and Lemma \ref{grad for all ui and ui conv to u impl grad for u} imply that
 $\widehat{g}_j$ is indeed a $p,q$-weak upper gradient for $w_j.$

 Moreover,
 \begin{equation} \label{gjk is ug for ujk}
 w_{j} \le v_{j}+ \sum_{i=j}^{\infty} |v_{i+1}-v_{i}| \mbox{ and }
 \widehat{g}_{j}\le \widetilde{g}_{j}+\sum_{i=j}^{k-1} |\widetilde{g}_{{i+1}}-\widetilde{g}_{i}|
 \end{equation}
 Thus
 \begin{eqnarray*}
 ||w_{j}||_{L^{p,q}(X, \mu)} &\le& ||v_{j}||_{L^{p,q}(X, \mu)} + \sum_{i=j}^{\infty} ||v_{{i+1}}-v_{i}||_{L^{p,q}(X, \mu)}
 \le ||v_{j}||_{L^{p,q}(X, \mu)}+ 2^{-j+1} \mbox{ and }\\
 ||\widehat{g}_{j}||_{L^{p,q}(X, \mu)} &\le& ||\widetilde{g}_{j}||_{L^{p,q}(X, \mu)} + \sum_{i=j}^{\infty} ||\widetilde{g}_{{i+1}}-\widetilde{g}_{i}||_{L^{p,q}(X, \mu)}
 \le ||\widetilde{g}_{j}||_{L^{p,q}(X, \mu)}+ 2^{-j+1},
 \end{eqnarray*}
 which implies that $w_j, \widehat{g}_j \in L^{p,q}(X, \mu).$  Thus $w_j \in {\mathcal{A}}(E)$ with $p,q$-weak upper gradient $\widehat{g}_j.$
 We notice that $0 \le g=\inf_{j \ge 1} \widehat{g}_{j}$ $\mu$-almost everywhere on $X$ and $0 \le u=\inf_{j \ge 1} w_j$ $\mu$-almost everywhere on $X.$ Since $w_1$ and $\widehat{g}_1$ are in $L^{p,q}(X,\mu),$ the absolute continuity of the $p,q$-norm (see
 Bennett-Sharpley \cite[Proposition I.3.6]{BS} and the discussion after Definition
 \ref{definition of abs continuous norm}) yields
 \begin{eqnarray} \label{wj conv to u and widehat gj conv to g}
 \lim_{j \rightarrow \infty} ||w_j-u||_{L^{p,q}(X,\mu)}=0 \mbox{ and } \lim_{j \rightarrow \infty} ||\widehat{g}_j-g||_{L^{p,q}(X,\mu)}=0.
 \end{eqnarray}
 By using (\ref{norm of (u,g) controlled by norm of L to some power}),
 (\ref{wj conv to u and widehat gj conv to g}), and
 Corollary \ref{Convergence of the quasinorm of the sequence} we see that
 $${\mathrm{Cap}}_{p,q}(E)^{q/p} \le
 \lim_{j \rightarrow \infty} ||w_{j}||_{N^{1,L^{p,q}}}^{q}= ||u||_{L^{p,q}(X, \mu)}^{q}+ ||g||_{L^{p,q}(X, \mu)}^{q} \le L^{q/p}+ \varepsilon.$$ By letting
 $\varepsilon \rightarrow 0,$ we get the converse inequality so (ii) is
 proved.

 \vskip 2mm

 (iii) We can assume without loss of generality that
 $$\sum_{i=1}^{\infty} {\mathrm{Cap}}_{p,q}(E_{i})^{q/p}<\infty.$$

 For $i=1,2, \ldots$ let $u_{i} \in {\mathcal{A}}(E_{i})$ with upper gradient $g_{i}$ such that
 $$0 \le u_{i} \le 1 \mbox{ and }
 ||u_{i}||_{N^{1,L^{p,q}}}^{q} < {\mathrm{Cap}}_{p,q}(E_{i})^{q/p}+ \varepsilon 2^{-i}.$$
 Let $u:=(\sum_{i=1}^{\infty} u_i^q)^{1/q}$ and $g:=(\sum_{i=1}^{\infty} g_i^q)^{1/q}.$ We notice that $u \ge 1$ on $E.$ By repeating the argument from the proof of Theorem \ref{basic properties of pq modulus} (iii), we see that $u,g \in L^{p,q}(X,\mu)$ and
 $$||u||_{L^{p,q}(X, \mu)}^q + ||g||_{L^{p,q}(X, \mu)}^q \le \sum_{i=1}^{\infty} \left(||u_{i}||_{L^{p,q}(X, \mu)}^q+||g_{i}||_{L^{p,q}(X, \mu)}^q\right) \le 2\varepsilon+\sum_{i=1}^{\infty} {\mathrm{Cap}}_{p,q}(E_{i})^{q/p}.$$
 We are done with the case $1 \le q \le p$ as soon as we show that $u \in {\mathcal{A}}(E)$ and that $g$ is a $p,q$-weak upper gradient for $u.$
 It follows easily via Corollary \ref{grad for root q of ulq plus u2q} and finite induction that $g$ is a $p,q$-weak upper gradient for $\widetilde{u}_n:=(\sum_{1 \le i \le n} u_i^q)^{1/q}$ for every $n \ge 1.$ Since $u(x)=\lim_{i \rightarrow \infty} \widetilde{u}_i(x)<\infty$ on $X \setminus F,$ where $F=\{ x \in X: u(x)=\infty \}$ it follows from  Lemma \ref{grad for all ui and ui conv to u impl grad for u} combined with the fact that $u \in L^{p,q}(X, \mu)$ that $g$ is in fact a $p,q$-weak upper gradient for $u.$ This finishes the proof for the case $1 \le q \le p.$

 \vskip 2mm

 (iv) We can assume without loss of generality that
 $$\sum_{i=1}^{\infty} {\mathrm{Cap}}_{p,q}(E_{i})<\infty.$$
 For $i=1, 2, \ldots$ let $u_{i} \in {\mathcal{A}}(E_{i})$ with upper gradients $g_{i}$ such that
 $$0 \le u_{i} \le 1 \mbox{ and }
 ||u_{i}||_{N^{1,L^{p,q}}}^{p} < {\mathrm{Cap}}_{p,q}(E_{i})+ \varepsilon 2^{-i}.$$

 Let $u:=\sup_{i \ge 1} u_{i}$ and $g:=\sup_{i \ge 1} g_i.$ We notice that $u=1$ on $E.$ Moreover, via Proposition \ref{subadd p le q} it follows that
 $u, g \in L^{p,q}(X,\mu)$ with
 $$||u||_{L^{p,q}(X, \mu)}^{p} + ||g||_{L^{p,q}(X, \mu)}^{p} \le \sum_{i=1}^{\infty} \left(||u_{i}||_{L^{p,q}(X, \mu)}^p+||g_{i}||_{L^{p,q}(X, \mu)}^p\right) \le 2\varepsilon+\sum_{i=1}^{\infty} {\mathrm{Cap}}_{p,q}(E_{i}).$$
 We are done with the case $p<q \le \infty$ as soon as we show that $u \in {\mathcal{A}}(E)$ and that $g$ is a $p,q$-weak upper gradient for $u.$
 Via Lemma \ref{gradient for max of two functions} and finite induction, it follows that $g$ is a $p,q$-weak upper gradient for $\widetilde{u}_n:=\max_{1 \le i \le n} u_i$ for every $n \ge 1.$ Since $u(x)=\lim_{i \rightarrow \infty} \widetilde{u}_i(x)$ pointwise on $X,$ it follows via Lemma \ref{grad for all ui and ui conv to u impl grad for u} that $g$ is in fact a $p,q$-weak upper gradient for $u.$ This finishes the proof for the case $p<q \le \infty.$

 \end{pf}

 \begin{Remark} \label{Remarks about the pq capacity} We make a few remarks.

 (i) Suppose $\mu$ is nonatomic and $1<q<\infty.$ By mimicking the proof of Theorem \ref{p,q Cap Thm} (ii) and working with the $(p,q)$-norm and the $(p,q)$-capacity, we can also show that
 $$ \lim_{i \rightarrow \infty} {\mathrm{Cap}}_{(p,q)}(E_i)= {\mathrm{Cap}}_{(p,q)}(E)$$
 whenever $E_{1} \subset E_{2} \subset \ldots \subset E=\bigcup_{i=1}^{\infty} E_{i} \subset X.$

 (ii) Moreover, if ${\mathrm{Cap}}_{p,q}$ is an outer capacity
 then it follows immediately that
 $$ \lim_{i \rightarrow \infty} {\mathrm{Cap}}_{p,q}(K_i)= {\mathrm{Cap}}_{p,q}(K)$$
 whenever $(K_i)_{i=1}^{\infty}$ is a decreasing sequence of compact sets whose intersection set is $K.$ We say that ${\mathrm{Cap}}_{p,q}$ is an outer capacity if for every $E \subset X$ we have
 $${\mathrm{Cap}}_{p,q}(E)= \inf \{ {\mathrm{Cap}}_{p,q}(U): E \subset U \subset X, \, U \mbox { open} \}.$$

 (iii) Any outer capacity satisfying properties (i) and (ii) of
 Theorem \ref{p,q Cap Thm} is called a Choquet capacity. (See Appendix II in Doob \cite{Doo}.)

 \end{Remark}

 We recall that if $A \subset X,$ then $\Gamma_{A}$ is the family of curves in $\Gamma_{\rm{rect}}$ that intersect $A$ and $\Gamma_{A}^{+}$ is the family of all curves in $\Gamma_{\rm{rect}}$ such that the Hausdorff one-dimensional measure ${\mathcal{H}}_{1}(|\gamma| \cap A)$ is positive. The following lemma will be useful later in this paper.

 \begin{Lemma} \label{cap pq F is zero implies mod pq gamma F is zero} If $F \subset X$ is such that ${\mathrm{Cap}}_{p,q}(F)=0,$ then ${\mathrm{Mod}}_{p,q}(\Gamma_F)=0.$

 \end{Lemma}

 \begin{pf} We follow Shanmugalingam \cite{Sha1}. We can assume without loss of generality that $q \neq p.$ Since ${\mathrm{Cap}}_{p,q}(F)=0,$ for each positive
 integer $i$ there exists a function $v_i \in {\mathcal{A}}(F)$ such that $0 \le v_i \le 1$ and such that $||v_i||_{N^{1,L^{(p,q)}}} \le 2^{-i}.$ Let $u_n:=\sum_{i=1}^{n} v_i.$
 Then $u_n \in N^{1,L^{(p,q)}}(X, \mu)$ for each $n,$ $u_n(x)$ is increasing for each $x \in X,$ and for every $m>n$ we have
 $$||u_n-u_m||_{N^{1,L^{(p,q)}}} \le \sum_{i=m+1}^{n} ||v_i||_{N^{1,L^{(p,q)}}} \le 2^{-m} \rightarrow 0, \mbox{ as } m \rightarrow \infty.$$

 Therefore the sequence $\{u_n\}_{n=1}^{\infty}$ is a Cauchy sequence in $N^{1,L^{(p,q)}}(X, \mu).$

 Since $\{u_n\}_{n=1}^{\infty}$ Cauchy in $N^{1,L^{(p,q)}}(X, \mu),$ it follows that it is Cauchy in $L^{p,q}(X, \mu).$ Hence by passing to a subsequence if necessary, there is a function $\widetilde{u}$ in $L^{p,q}(X, \mu)$ to which the subsequence converges both pointwise $\mu$-almost everywhere and in the $L^{(p,q)}$ norm. By choosing a further subsequence, again denoted by $\{u_i\}_{i=1}^{\infty}$ for simplicity, we can assume that
 $$||u_i-\widetilde{u}||_{L^{(p,q)}(X, \mu)} +||g_{i, i+1}||_{L^{(p,q)}(X, \mu)} \le 2^{-2i},$$ where $g_{i,j}$ is an upper gradient of $u_i-u_j$ for $i<j.$ If $g_1$ is an upper gradient of $u_1,$ then $u_2=u_1+(u_2-u_1)$ has an upper gradient $g_2=g_1+g_{12}.$ In general,
 $$u_i=u_{1}+ \sum_{k=1}^{i-1} (u_{k+1}-u_{k}) $$ has an upper gradient
 $$g_i=g_{1}+ \sum_{k=1}^{i-1} g_{k, k+1}$$ for every $i \ge 2.$
 For $j<i$ we have
 $$||g_{i}-g_{j}||_{L^{(p,q)}(X, \mu)} \le \sum_{k=j}^{i-1} ||g_{k, k+1}||_{L^{(p,q)}(X, \mu)} \le \sum_{k=j}^{i-1} 2^{-2k} \le 2^{1-2j} \rightarrow 0 \mbox{ as } j \rightarrow \infty.$$

 Therefore $\{g_i\}_{i=1}^{\infty}$ is also a Cauchy sequence in $L^{(p,q)}(X, \mu),$ and hence converges in the $L^{(p,q)}$ norm to a nonnegative Borel function $g.$
 Moreover, we have $$||g_{j}-g||_{L^{(p,q)}(X, \mu)} \le 2^{1-2j}$$
 for every $j \ge 1.$

 We define $u$ by $u(x)=\lim_{i \rightarrow \infty} u_{i}(x).$ Since $u_{i} \rightarrow \widetilde{u}$ $\mu$-almost everywhere, it follows that $u=\widetilde{u}$ $\mu$-almost everywhere and thus $u \in L^{p,q}(X, \mu).$ Let $$E= \{ x \in X: \lim_{i \rightarrow \infty} u_{i}(x)=\infty \}.$$ The function $u$ is well defined outside of $E.$ In order for the function $u$ to be in the space $N^{1, L^{p,q}}(X, \mu),$ the function $u$ has to be defined on almost all paths by Proposition \ref{N1pq is a subset of ACCpq}. To this end it is shown that the $p,q$-modulus of the family $\Gamma_{E}$ is zero. Let $\Gamma_{1}$ be the collection of all paths from $\Gamma_{\rm{rect}}$ such that $\int_{\gamma} g=\infty.$ Then we have via Theorem
 \ref{char of families of curves of pq modulus zero} that ${\mathrm{Mod}}_{p,q}(\Gamma_1)=0$ since $g \in L^{p,q}(X, \mu).$

 Let $\Gamma_2$ be the family of all curves from $\Gamma_{\rm{rect}}$ such that $\limsup_{j \rightarrow \infty} \int_{\gamma} |g_j-g|>0.$ Since $||g_j-g||_{L^{p,q}}(X,\mu) \le 2^{1-2j}$ for all $j \ge 1,$ it follows via Theorem
 \ref{uk conv to u in Lpq implies uk conv to u on curves} that ${\mathrm{Mod}}_{p,q}(\Gamma_2)=0.$

 Since $u \in L^{p,q}(X, \mu)$ and $E=\{x \in X: u(x)=\infty \},$ it follows that $\mu(E)=0$ and thus ${\mathrm{Mod}}_{\Gamma_{E}^{+}}=0.$ Therefore ${\mathrm{Mod}}_{p,q}(\Gamma_1 \cup \Gamma_2 \cup \Gamma_{E}^{+})=0.$ For any path $\gamma$ in the family $\Gamma_{\rm{rect}} \setminus (\Gamma_1 \cup \Gamma_2 \cup \Gamma_{E}^{+}),$ by the fact that $\gamma$ is not in $\Gamma_{E}^{+},$ there exists a point in $|\gamma| \setminus E.$ For any point $x$ in $\gamma,$ since $g_i$ is an upper gradient of $u_i,$ it follow that
 $$u_{i}(x)-u_{i}(y) \le |u_{i}(x)-u_{i}(y)| \le \int_{\gamma} g_i.$$
 Therefore $$ u_{i}(x) \le u_{i}(y)+ \int_{\gamma} g_i.$$
 Taking limits on both sides and using the fact that $\gamma$ is not in $\Gamma_1 \cup \Gamma_2,$ it follows that
 $$\lim_{i \rightarrow \infty} u_{i}(x) \le \lim_{i \rightarrow \infty} u_{i}(y)+ \int_{\gamma} g = u(y)+ \int_{\gamma} g < \infty,$$ and therefore $x$ is not in $E.$ Thus $\Gamma_{E} \subset \Gamma_1 \cup \Gamma_2 \cup \Gamma_{E}^{+}$ and ${\mathrm{Mod}}_{p,q}(\Gamma_{E})=0.$
 Therefore $g$ is a $p,q$-weak upper gradient of $u,$ and hence $u \in N^{1, L^{p,q}}(X, \mu).$ For each $x$ not in $E$ we can write $u(x)=\lim_{i \rightarrow \infty} u_{i}(x)<\infty.$ If $F \setminus E$ is nonempty, then
 $$u|_{F \setminus E} \ge u_n|_{F \setminus E}= \sum_{i=1}^{n} v_i|_{F \setminus E}=n$$
 for arbitrarily large $n,$ yielding that $u|_{F \setminus E}=\infty.$ But this impossible, since $x$ is not in the set $E.$ Therefore $F \subset E,$ and hence
 $\Gamma_{F} \subset \Gamma_{E}.$ This finishes the proof of the lemma.

 \end{pf}

 Next we prove that $(N^{1, L^{p,q}}(X, \mu), ||\cdot||_{N^{1, L^{(p,q)}}})$ is a Banach space.

 \begin{Theorem} \label{Newtonian pq space is Banach} Suppose $1<p<\infty$ and $1 \le q \le \infty.$ Then $(N^{1, L^{p,q}}(X,\mu), ||\cdot||_{N^{1, L^{(p,q)}}})$ is a Banach space.

 \end{Theorem}

 \begin{pf} We follow Shanmugalingam \cite{Sha1}. We can assume without loss of generality that $q \neq p.$ Let $\{u_i\}_{i=1}^{\infty}$
 be a Cauchy sequence in $N^{1, L^{p,q}}(X,\mu).$ To show that this sequence is convergent in $N^{1, L^{p,q}}(X,\mu),$ it suffices to show that some subsequence is convergent in $N^{1, L^{p,q}}(X,\mu).$ Passing to a further subsequence if necessary, it can be assumed that
 $$||u_{i+1}-u_{i}||_{L^{(p,q)}(X,\mu)} +||g_{i, i+1}||_{L^{(p,q)}(X,\mu)} \le 2^{-2i},$$
 where $g_{i,j}$ is an upper gradient of $u_i-u_j$ for $i<j.$ Let
 $$E_{j}= \{ x \in X: |u_{j+1}(x)-u_{j}(x)| \ge 2^{-j} \}.$$
 Then $2^{j}|u_{j+1}-u_{j}| \in {\mathcal{A}}(E_{j})$ and hence
 $${\mathrm{Cap}}_{p,q}(E_{j})^{1/p} \le 2^{j} ||u_{j+1}-u_{j}||_{N^{1, L^{p,q}}} \le 2^{-j}.$$
 Let $F_{j}=\cup_{k=j}^{\infty} E_{k}.$ Then
 $${\mathrm{Cap}}_{p,q}(E_{j})^{1/p} \le \sum_{k=j}^{\infty} {\mathrm{Cap}}_{p,q}(E_{k})^{1/p} \le 2^{1-j}.$$
 Let $F= \cap_{j=1}^{\infty} F_{j}.$ We notice that ${\mathrm{Cap}}_{p,q}(F)=0.$
 If $x$ is a point in $X \setminus F,$ there exists $j \ge 1$ such that $x$ is not in
 $F_{j}=\cup_{k=j}^{\infty} E_{k}.$ Hence for all $k \ge j,$ $x$ is not in $E_{k}.$ Thus
 $|u_{k+1}(x)-u_{k}(x)| \le 2^{-k}$ for all $k \ge j.$ Therefore whenever $l \ge k \ge j$ we have that
 $$|u_{k}(x)-u_{l}(x)| \le 2^{1-k}.$$
 Thus the sequence $\{u_{k}(x)\}_{k=1}^{\infty}$ is Cauchy for every $x \in X \setminus F.$ For every $x \in X \setminus F,$ let $u(x)=\lim_{i \rightarrow \infty} u_{i}(x).$
 For $k<m,$
 $$u_{m}=u_{k}+ \sum_{n=k}^{m-1}(u_{n+1}-u_{n}).$$
 Therefore for each $x$ in $X \setminus F,$
 \begin{equation} \label{u minus uk equals sum unplus1 minus un}
 u(x)=u_{k}(x)+ \sum_{n=k}^{\infty}(u_{n+1}(x)-u_{n}(x)).
 \end{equation}

 Noting by Lemma \ref{cap pq F is zero implies mod pq gamma F is zero} that ${\mathrm{Mod}}_{p,q}(\Gamma_{F})=0$ and that for each path $\gamma$ in $\Gamma_{\rm{rect}} \setminus \Gamma_{F}$ equation
 (\ref{u minus uk equals sum unplus1 minus un}) holds pointwise on $|\gamma|,$ we conclude that $\sum_{n=k}^{\infty} g_{n, n+1}$ is a $p,q$-weak upper gradient of $u-u_k.$
 Therefore
 \begin{eqnarray*}
 ||u-u_{k}||_{N^{1, L^{(p,q)}}} &\le& ||u-u_{k}||_{L^{(p,q)}(X, \mu)}+ \sum_{n=k}^{\infty} ||g_{n, n+1}||_{L^{(p,q)}(X, \mu)} \\
 &\le& ||u-u_k||_{L^{(p,q)}(X,\mu)}+ \sum_{n=k}^{\infty} 2^{-2n} \\
 &\le& ||u-u_k||_{L^{(p,q)}(X,\mu)}+2^{1-2k} \rightarrow 0 \mbox{ as } k \rightarrow \infty.
 \end{eqnarray*}
 Therefore the subsequence converges in the norm of $N^{1,L^{p,q}}(X,\mu)$ to $u.$
 This completes the proof of the theorem.

 \end{pf}

 \section{Density of Lipschitz functions in $N^{1, L^{p,q}}(X, \mu)$}

 \subsection{Poincar\'{e} inequality}

 Now we define the weak $(1, L^{p,q})$-\textit{Poincar\'{e} inequality}. Podbrdsky in \cite{Pod} introduced a stronger Poincar\'{e} inequality in the case of Banach-valued Newtonian Lorentz spaces.

 \begin{Definition} \label{definition of 1pq Poincare inequality}
 The space $(X,d,\mu)$ is said to support a \textit{weak $(1, L^{p,q})$-Poincar\'{e} inequality} if there exist constants $C>0$ and $\sigma \ge 1$ such that for all balls $B$ with radius $r,$ all $\mu$-measurable functions $u$ on $X$ and all upper gradients $g$ of $u$ we have
 \begin{equation} \label{1pq Poincare inequality}
 \frac{1}{\mu(B)} \int_{B} |u-u_{B}| \, d\mu \le
 Cr \frac{1}{\mu(\sigma B)^{1/p}} ||g \chi_{\sigma B}||_{L^{p,q}(X, \mu)}.
 \end{equation}
 Here $$u_{B}=\frac{1}{\mu(B)} \int_{B} u(x) \, d\mu(x)$$
 whenever $u$ is a locally $\mu$-integrable function on $X.$
 \end{Definition}

 In the above definition we can equivalently assume via Lemma \ref{pq weak upper gradient approximated by upper gradients} and Corollary \ref{Convergence of the quasinorm of the sequence} that $g$ is a $p,q$-weak upper gradient of $u.$ When $p=q$ we have the weak $(1,p)$-Poincar\'{e} inequality. For more about the Poincar\'{e} inequality in the case $p=q$ see Haj{\l}asz-Koskela \cite{HaK} and Heinonen-Koskela \cite{HeK}.

 A measure $\mu$ is said to be \textit{doubling} if there exists a constant $C\ge 1$ such that
 $$\mu(2B) \le C \mu(B)$$
 for every ball $B=B(x,r)$ in X. A metric measure space $(X,d,\mu)$ is called \textit{doubling}
 if the measure $\mu$ is doubling. Under the assumption that the measure $\mu$ is doubling,
 it is known that $(X,d,\mu)$ is proper (that is, closed bounded subsets of $X$ are compact) if and only if it is complete.

 Now we prove that if $1 \le q \le p,$ the measure $\mu$ is doubling, and the space $(X,d,\mu)$ carries a weak $(1, L^{p,q})$-Poincar\'{e} inequality, the Lipschitz functions are dense in $N^{1, L^{p,q}}(X, \mu).$

 In order to prove this we need a few definitions and lemmas.

 \begin{Definition} Suppose $(X, d)$ is a metric space that carries a doubling measure $\mu.$ For $1 <p <\infty$ and
 $1 \le q \le \infty$ we define the noncentered maximal function operator
 by $$M_{p,q}u(x)=\sup_{B \ni x} \frac{||u \chi_{B}||_{L^{p,q}(X, \mu)}}{\mu(B)^{1/p}},$$
 where $u \in L^{p,q}(X,\mu).$

 \end{Definition}

 \begin{Lemma} \label{Mpq maps Lpq to weak Lp for qlep} Suppose $(X, d)$ is a metric space that carries a doubling measure $\mu.$ If $1 \le q \le p,$ then $M_{p,q}$ maps $L^{p,q}(X,\mu)$ to $L^{p,\infty}(X,\mu)$ boundedly and moreover, $$\lim_{\lambda \rightarrow \infty} \lambda^{p} \mu(\{ x \in X: M_{p,q}u(x)>\lambda \})=0.$$

 \end{Lemma}

 \begin{pf} We can assume without loss of generality that $1 \le q<p.$ For every $R>0$ let $M_{p,q}^{R}$ be the restricted maximal function operator defined on $L^{p,q}(X, \mu)$ by
 $$M_{p,q}^{R}u(x)=\sup_{B \ni x, \mbox{ diam}(B) \le R} \frac{||u \chi_{B}||_{L^{p,q}(X, \mu)}}{\mu(B)^{1/p}}.$$

 Denote $G_{\lambda}=\{ x \in X: M_{p,q}u(x)>\lambda \}$ and
 $G_{\lambda}^{R}=\{ x \in X: M_{p,q}^{R}u(x)>\lambda \}.$ It is easy to see that $G_{\lambda}^{R_1} \subset G_{\lambda}^{R_2}$ if $0<R_1<R_2<\infty$ and
 $G_{\lambda}^{R} \rightarrow G_{\lambda}$ as $R \rightarrow \infty.$

 Fix $R>0.$ For every $x \in G_{\lambda}^{R},$ $\lambda>0,$ there exists a ball $B(y_{x}, r_{x})$ with diameter at most $R$ such that $x \in B(y_x, r_x)$ and such that
 $$||u \chi_{B(y_x, r_x)}||_{L^{p,q}(X, \mu)}^p > \lambda^p \mu(B(y_x, r_x)).$$
 We notice that $B(y_x, r_x) \subset G_{\lambda}^{R}.$ The set $G_{\lambda}^{R}$ is covered by such balls and by Theorem 1.2 in Heinonen \cite{Hei} it follows that there exists a countable disjoint subcollection $\{ B(x_i, r_i) \}_{i=1}^{\infty}$ such that
 the collection $\{ B(x_i, 5r_i) \}_{i=1}^{\infty}$ covers $G_{\lambda}^{R}.$ Hence
 \begin{eqnarray*}
 \mu(G_{\lambda}^{R}) &\le& \sum_{i=1}^{\infty} \mu(B(x_i, 5r_i)) \le C \left(\sum_{i=1}^{\infty} \mu(B(x_i, r_i))\right)\\
 &\le& \frac{C}{\lambda^{p}} \left( \sum_{i=1}^{\infty} ||u \chi_{B(x_i, r_i)}||_{L^{p,q}(X, \mu)}^p \right) \le \frac{C}{\lambda^{p}} ||u \chi_{G_{\lambda}^{R}}||_{L^{p,q}(X, \mu)}^{p}.
 \end{eqnarray*}
 The last inequality in the sequence was obtained by applying Proposition \ref{superadd q le p}. (See also Chung-Hunt-Kurtz \cite[Lemma 2.5]{CHK}.)

 Thus $$\mu(G_{\lambda}^{R}) \le \frac{C}{\lambda^{p}} ||u \chi_{G_{\lambda}^{R}}||_{L^{p,q}(X, \mu)}^{p} \le \frac{C}{\lambda^{p}} ||u \chi_{G_{\lambda}}||_{L^{p,q}(X, \mu)}^{p}$$
 for every $R>0.$ Since $G_{\lambda}=\bigcup_{R>0} G_{\lambda}^{R},$ we obtain (by taking the limit as $R \rightarrow \infty$)
 $$\mu(G_{\lambda}) \le \frac{C}{\lambda^{p}} ||u \chi_{G_{\lambda}}||_{L^{p,q}(X, \mu)}^{p}.$$
 The absolute continuity of the $p,q$-norm (see the discussion after Definition \ref{definition of abs continuous norm}), the $p,q$-integrability of $u$ and the fact that $G_{\lambda} \rightarrow \emptyset$ $\mu$-almost everywhere as $\lambda \rightarrow \infty$ yield the desired conclusion.

 \end{pf}

 \begin{Question} Is Lemma \ref{Mpq maps Lpq to weak Lp for qlep} true when $p<q< \infty$ ?

 \end{Question}

 The following proposition is necessary in the sequel.

 \begin{Proposition} \label{sequence of truncations conv in N1pq q finite} Suppose $1<p<\infty$ and $1 \le q<\infty.$ If $u$ is a nonnegative function in $N^{1,L^{p,q}}(X,\mu),$ then the sequence of functions $u_k=\min(u,k),$
 $k \in \Bbb{N},$ converges in the norm of $N^{1, L^{p,q}}(X, \mu)$ to $u$ as $k \rightarrow \infty.$

 \end{Proposition}

 \begin{pf} We notice (see Lemma \ref{gradient truncation above}) that $u_k \in L^{p,q}(X, \mu).$ That lemma also yields easily $u_k \in N^{1, L^{p,q}}(X, \mu)$ and moreover $||u_k||_{N^{1,L^{p,q}}} \le ||u||_{N^{1,L^{p,q}}}$ for all $k \ge 1.$

 Let $E_{k}=\{ x \in X: u(x)>k \}.$ Since $\mu$ is a Borel regular measure, there exists an open set $O_k$ such that $E_k \subset O_k$ and $\mu(O_k) \le \mu(E_k)+ 2^{-k}.$
 In fact the sequence $(O_k)_{k=1}^{\infty}$ can be chosen such that $O_{k+1} \subset O_{k}$ for all $k \ge 1.$ Since $\mu(E_k) \le \frac{C(p,q)}{k^p} ||u||_{L^{p,q}(X, \mu)}^p,$ it follows that
 $$\mu(O_k) \le \mu(E_k) + 2^{-k} \le \frac{C(p,q)}{k^p} ||u||_{L^{p,q}(X, \mu)}^p+2^{-k}.$$
 Thus $\lim_{k \rightarrow \infty} \mu(O_k)=0.$ We notice that $u=u_k$ on $X \setminus O_k.$ Thus $2 g \chi_{O_k}$ is a $p,q$-weak upper gradient of $u-u_k$ whenever $g$
 is an upper gradient for $u$ and $u_k.$ See Lemma
 \ref{gradient trimming on closed sets}. The fact that $O_{k} \rightarrow \emptyset$ $\mu$-almost everywhere and the absolute continuity of the $(p,q)$-norm yield
 $$\limsup_{k \rightarrow \infty} ||u-u_{k}||_{N^{1,L^{(p,q)}}} \le 2 \limsup_{k \rightarrow \infty} \left(||u \chi_{O_k}||_{L^{(p,q)}(X, \mu)} + ||g \chi_{O_k}||_{L^{(p,q)}(X, \mu)}\right)=0. $$

 \end{pf}

 \begin{Counterexample} \label{counterexample conv truncation for N1pinfinity} For $q=\infty$ Proposition \ref{sequence of truncations conv in N1pq q finite} is not true. Indeed, let $n \ge 2$ be an integer and let $1<p \le n$ be fixed. Let $X=B(0,1) \setminus \{0\} \subset {\Bbb{R}}^{n},$ endowed with the Euclidean metric and the Lebesgue measure.

 Suppose first that $1<p<n.$ Let $u_p$ and $g_p$ be defined on $X$ by
 $$ u_{p}(x)= |x|^{1-\frac{n}{p}}-1 \mbox{ and } g_{p}(x)=\left(\frac{n}{p}-1 \right) |x|^{-\frac{n}{p}}.$$
 It is easy to see that $u_{p}, g_{p} \in L^{p, \infty}(X, m_n).$ Moreover, (see for instance Haj{\l}asz \cite[Proposition 6.4]{Haj2}) $g_{p}$ is the minimal upper gradient for $u_{p}.$ Thus $u_{p} \in N^{1,L^{p, \infty}}(X, m_n).$ For every integer $k \ge 1$ we define $u_{p,k}$ and $g_{p,k}$ on $X$ by
  $$u_{p,k}(x)=\left\{ \begin{array}{ll}
k & \mbox{if } 0<|x|\le (k+1)^{\frac{p}{p-n}}, \\
|x|^{1-\frac{n}{p}}-1 & \mbox{if } (k+1)^{\frac{p}{p-n}} < |x| < 1
\end{array}
\right.
$$

and  $$g_{p,k}(x)=\left\{ \begin{array}{ll}
\left(\frac{n}{p}-1\right)|x|^{-\frac{n}{p}} & \mbox{if } 0<|x|<(k+1)^{\frac{p}{p-n}}\\
0 & \mbox{if } (k+1)^{\frac{p}{p-n}} \le |x|<1.\\
\end{array}
\right.
$$
We notice that $u_{p,k} \in N^{1, L^{p,\infty}}(X, m_n)$ for all $k \ge 1.$ Moreover, via \cite[Proposition 6.4]{Haj2} and Lemma \ref{gradient trimming on closed sets} we see that $g_{p,k}$ is the minimal upper gradient for $u_{p}-u_{p,k}$ for every $k \ge 1.$ Since $g_{p,k} \searrow 0$ on $X$ as $k \rightarrow \infty$ and $||g_{p,k}||_{L^{p,\infty}(X, m_n)}=||g_{p}||_{L^{p,\infty}(X, m_n)}=C(n,p)>0$ for all $k \ge 1,$ it follows that
$u_{p,k}$ does \textit{not} converge to $u_{p}$ in $N^{1, L^{p,\infty}}(X, m_n)$ as $k \rightarrow \infty.$

Suppose now that $p=n.$ Let $u_n$ and $g_n$ be defined on $X$ by
 $$ u_{n}(x)= \ln \frac{1}{|x|} \mbox{ and } g_{n}(x)=\frac{1}{|x|}.$$
 It is easy to see that $u_{n}, g_{n} \in L^{p, \infty}(X, m_n).$ Moreover, (see for instance Haj{\l}asz \cite[Proposition 6.4]{Haj2}) $g_{n}$ is the minimal upper gradient for $u_{n}.$ Thus $u_{n} \in N^{1,L^{n, \infty}}(X, m_n).$ For every integer $k \ge 1$ we define $u_{n,k}$ and $g_{n,k}$ on $X$ by
  $$u_{n,k}(x)=\left\{ \begin{array}{ll}
k & \mbox{if } 0<|x|\le e^{-k}, \\
\ln \frac{1}{|x|} & \mbox{if } e^{-k} < |x| < 1
\end{array}
\right.
$$

and  $$g_{n,k}(x)=\left\{ \begin{array}{ll}
\frac{1}{|x|} & \mbox{if } 0<|x|<e^{-k}\\
0 & \mbox{if } e^{-k} \le |x|<1.\\
\end{array}
\right.
$$
We notice that $u_{n,k} \in N^{1, L^{n,\infty}}(X, m_n)$ for all $k \ge 1.$ Moreover, via \cite[Proposition 6.4]{Haj2} and Lemma \ref{gradient trimming on closed sets} we see that $g_{n,k}$ is the minimal upper gradient for $u_{n}-u_{n,k}$ for every $k \ge 1.$ Since $g_{n,k} \searrow 0$ on $X$ as $k \rightarrow \infty$ and $||g_{n,k}||_{L^{p,\infty}(X, m_n)}=||g_{n}||_{L^{n,\infty}(X, m_n)}=C(n)>0$ for all $k \ge 1,$ it follows that
$u_{n,k}$ does \textit{not} converge to $u_{n}$ in $N^{1, L^{n,\infty}}(X, m_n)$ as $k \rightarrow \infty.$

 \end{Counterexample}

 The following lemma will be used in the paper.

 \begin{Lemma} \label{product rule for N1pq} Let $f_1 \in N^{1, L^{p,q}}(X, \mu)$
 be a bounded Borel function with $p,q$-weak upper gradient $g_1 \in L^{p,q}(X, \mu)$ and let $f_2$ be a bounded Borel function with $p,q$-weak upper gradient $g_2 \in L^{p,q}(X,\mu).$
 Then $f_3:=f_1 f_2 \in N^{1, L^{p,q}}(X,\mu)$ and $g_3:=|f_1|g_2+|f_2|g_1 $ is a $p,q$-weak upper gradient of $f_3.$

 \end{Lemma}

 \begin{pf} It is easy to see that $f_3$ and $g_3$ are in $L^{p,q}(X,\mu).$
 Let $\Gamma_{0} \subset \Gamma_{\rm{rect}}$ be the family of curves on which $\int_{\gamma} (g_1+g_2) =\infty.$ Then it follows via Theorem
 \ref{char of families of curves of pq modulus zero} that ${\mathrm{Mod}}_{p,q}(\Gamma_{0})=0$ because $g_1+g_2 \in L^{p,q}(X,\mu).$

 Let $\Gamma_{1,i} \subset \Gamma_{\rm{rect}}, i=1,2$ be the family of curves for which
 $$|f_{i}(\gamma(0))-f_{i}(\gamma(\ell(\gamma)))| \le \int_{\gamma} g_{i} $$
 is \textit{not} satisfied. Then ${\mathrm{Mod}}_{\Gamma_{1,i}}=0, i=1,2.$
 Let $\Gamma_{1} \subset \Gamma_{\rm{rect}}$ be the family of curves that have a subcurve in $\Gamma_{1,1} \cup \Gamma_{1,2}.$ Then $F(\Gamma_{1,1} \cup \Gamma_{1,2}) \subset F(\Gamma_{1})$ and thus ${\mathrm{Mod}}_{p,q}(\Gamma_{1}) \le {\mathrm{Mod}}_{p,q}(\Gamma_{1,1} \cup \Gamma_{1,2})=0.$
 We notice immediately that ${\mathrm{Mod}}_{p,q}(\Gamma_{0} \cup \Gamma_{1})=0.$

 Fix $\varepsilon>0.$ By using the argument from Lemma 1.7 in Cheeger \cite{Che}, we see that $$|f_3(\gamma(0))-f_3(\gamma(\ell(\gamma)))|\le \int_{0}^{\ell(\gamma)} (|f_1(\gamma(s))|+\varepsilon)g_2(\gamma(s))+ (|f_2(\gamma(s))|+\varepsilon)g_1(\gamma(s)) \, ds$$
 for every curve $\gamma$ in $\Gamma_{\rm{rect}} \setminus (\Gamma_{0} \cup \Gamma_{1}).$ Letting $\varepsilon \rightarrow 0$ we obtain the desired claim.

 \end{pf}

 Fix $x_{0} \in X.$ For each integer $j>1$ we consider the function
 $$\eta_j(x)=\left\{ \begin{array}{ll}
1 & \mbox{if } d(x_0,x)\le j-1, \\
j-d(x_0,x) & \mbox{if } j-1 < d(x_0,x) \le j,\\
0 & \mbox{if } d(x_0,x)> j.
\end{array}
\right.
$$

 \begin{Lemma} \label{u etaj converges to in N1pq when q is finite} Suppose $1 \le q<\infty.$ Let $u$ be a bounded function in $N^{1, L^{p,q}}(X, \mu).$ Then the function $v_j=u \eta_j$ is also in $N^{1, L^{p,q}}(X, \mu).$ Furthermore, the sequence $v_j$ converges to $u$ in $N^{1, L^{p,q}}(X, \mu).$

 \end{Lemma}

 \begin{pf} We can assume without loss of generality that $u \ge 0.$ Let $g \in L^{p,q}(X,\mu)$ be an upper gradient for $u.$ It is easy to see by invoking Lemma \ref{gradient trimming on closed sets} that $h_j:=\chi_{B(x_0, j) \setminus{\overline{B}(x_0, j-1)}}$ is a $p,q$-weak upper gradient for $\eta_j$ and for $1-\eta_{j}.$ By using Lemma \ref{product rule for N1pq}, we see that $v_j \in N^{1,L^{p,q}}(X,\mu)$ and
 that $g_j:=u h_j+ g \eta_j$ is a $p,q$-weak upper gradient for $v_j.$  By using Lemma \ref{product rule for N1pq}  we notice that $\widetilde{h}_{j}:=u h_{j}+ g (1-\eta_{j})$ is a $p,q$-weak upper gradient for $u-v_j.$ We have in fact
 \begin{equation} \label{upper est for u minus vj and for grad of u minus vj}
 0 \le u-v_j \le u \chi_{X \setminus B(x_0, j-1)} \mbox{ and } \widetilde{h}_{j} \le (u+g) \chi_{X \setminus B(x_0, j-1)}.
 \end{equation}
 for every $j>1.$ The absolute continuity of the $(p,q)$-norm when $1 \le q<\infty$ (see the discussion after Definition \ref{definition of abs continuous norm}) together with the $p,q$-integrability of $u, g$ and (\ref{upper est for u minus vj and for grad of u minus vj}) yield the desired claim.

  \end{pf}

 Now we prove the density of the Lipschitz functions in $N^{1, L^{p,q}}(X, \mu)$ when
 $1 \le q < p.$ The case $q=p$ was proved by Shanmugalingam. (See \cite{Sha1} and \cite{Sha2}.)

 \begin{Theorem} \label{Lip is dense in N1pq when q le p} Let $1\le q \le p<\infty.$ Suppose that $(X,d,\mu)$ is a doubling metric measure space that carries a weak $(1, L^{p,q})$-Poincar\'{e} inequality. Then the Lipschitz functions are dense in $N^{1, L^{p,q}}(X,\mu).$

 \end{Theorem}

 \begin{pf} We can consider only the case $1 \le q<p$ because the case $q=p$ was proved by Shanmugalingam in \cite{Sha1} and \cite{Sha2}. We can assume without loss of generality that $u$ is nonnegative. Moreover, via Lemmas \ref{sequence of truncations conv in N1pq q finite} and \ref{product rule for N1pq} we can assume without loss of generality that $u$ is bounded and has compact support in $X.$ Choose $M>0$ such that $0 \le u \le M$ on $X.$ Let $g \in L^{p,q}(X,\mu)$ be a $p,q$-weak upper gradient for $u.$ Let $\sigma \ge 1$ be the constant from the weak $(1, L^{p,q})$-Poincar\'{e} inequality.

 Let $G_{\lambda}:=\{ x \in X: M_{p,q}g(x)>\lambda \}.$ If $x$ is a point in the closed set $X \setminus G_{\lambda},$ then for all $r>0$ one has that
 $$\frac{1}{\mu(B(x,r))} \int_{B(x,r)} |u-u_{B(x,r)}| \, d\mu \le C r \frac{||g \chi_{B(x, \sigma r)}||_{L^{p,q}(X,\mu)}}{\mu(B(x, \sigma r))^{1/p}} \le C r M_{p,q}g(x) \le C \lambda r.$$
 Hence for $s \in [r/2, r]$ one has that
 \begin{eqnarray*}
 |u_{B(x,s)}-u_{B(x,r)}| &\le& \frac{1}{\mu(B(x,s))} \int_{B(x,s)} |u-u_{B(x,r)}| \, d\mu\\
 &\le& \frac{\mu(B(x,r))}{\mu(B(x,s))} \cdot \frac{1}{\mu(B(x,r))} \int_{B(x,r)} |u-u_{B(x,r)}| \, d\mu \le C \lambda r
 \end{eqnarray*}
 whenever $x$ is in $X \setminus G_{\lambda}.$
 For a fixed $s \in (0, r/2)$ there exists an integer $k \ge 1$ such that $2^{-k}r \le 2s<2^{-k+1}r.$
 Then
 \begin{eqnarray*}
 |u_{B(x,s)}-u_{B(x,r)}| &\le& |u_{B(x,s)}-u_{B(x,2^{-k}r)}|+     \\
 &+& \sum_{i=0}^{k-1} |u_{B(x, 2^{-i-1}r)}-u_{B(x,2^{-i}r)}| \le C \lambda \left(\sum_{i=0}^{k} 2^{-i}r\right) \le C \lambda r.
 \end{eqnarray*}
 For any sequence $r_i \searrow 0$ we notice that $(u_{B(x,r_i)})_{i=1}^{\infty}$ is a Cauchy sequence for every point $x$ in $X \setminus G_{\lambda}.$
 Thus on $X \setminus G_{\lambda}$ we can define the function
 $$u_{\lambda}(x):=\lim_{r \rightarrow 0} u_{B(x,r)}.$$
 We notice that $u_{\lambda}(x)=u(x)$ for every Lebesgue point $x$ in $X \setminus G_{\lambda}.$

 For every $x, y$ in $X \setminus G_{\lambda}$ we consider the chain of balls $\{B_i\}_{i=-\infty}^{\infty},$ where
 $$B_{i}=B(x, 2^{1+i}d(x,y)), i \le 0 \mbox{ and } B_{i}=B(y, 2^{1-i}d(x,y)), i>0.$$
 For every two such points $x$ and $y$ we have that they are Lebesgue points for
 $u_{\lambda}$ by construction and hence
 $$|u_{\lambda}(x)-u_{\lambda}(y)| \le \sum_{i=-\infty}^{\infty}|u_{B_{i}}-u_{B_{i+1}}| \le C \lambda d(x,y),$$
 where $C$ depends only on the data on $X.$ Thus $u_{\lambda}$ is $C \lambda$-Lipschitz on $X \setminus G_{\lambda}.$ By construction it follows that $0 \le u_{\lambda} \le M$
 on $X \setminus G_{\lambda}.$ Extend $u_{\lambda}$ as a $C \lambda$-Lipschitz function on $X$ (see McShane \cite{MS}) and denote the extension by $v_{\lambda}.$ Then $v_{\lambda} \ge 0$ on $X$ since $u_{\lambda} \ge 0$ on $X \setminus G_{\lambda}.$ Let $w_{\lambda}:=\min(v_{\lambda}, M).$ We notice that $w_{\lambda}$ is a nonnegative $C\lambda$-Lipschitz function on $X$ since $v_{\lambda}$ is. Moreover, $w_{\lambda}=v_{\lambda}=u_{\lambda}$ on $X \setminus G_{\lambda}$ whenever $\lambda>M.$

 Since $u=w_{\lambda}$ $\mu$-almost everywhere on $X \setminus G_{\lambda}$ whenever $\lambda>M$ we have
 $$||u-w_{\lambda}||_{L^{p,q}(X, \mu)}=||(u-w_{\lambda}) \chi_{G_{\lambda}}||_{L^{p,q}(X,\mu)} \le ||u \chi_{G_{\lambda}}||_{L^{p,q}(X,\mu)} + C(p,q) \lambda \mu(G_{\lambda})^{1/p}$$
 whenever $\lambda>M.$
 The absolute continuity of the $p,q$-norm when $1 \le q \le p$ together with Lemma \ref{Mpq maps Lpq to weak Lp for qlep} imply that $$ \lim_{\lambda \rightarrow \infty} ||u-w_{\lambda}||_{L^{p,q}(X, \mu)}=0.$$

 Since $u-w_{\lambda}=0$ $\mu$-almost everywhere on the closed set $G_{\lambda},$ it follows via Lemma \ref{gradient trimming on closed sets} that $(C \lambda+g) \chi_{G_{\lambda}}$ is a $p,q$-weak upper gradient for $u-w_{\lambda}.$ By using the absolute continuity of the $p,q$-norm when $1 \le q \le p$ together with Lemma \ref{Mpq maps Lpq to weak Lp for qlep}, we see that $$ \lim_{\lambda \rightarrow \infty} ||(C \lambda+g)\chi_{G_{\lambda}}||_{L^{p,q}(X, \mu)}=0.$$

 This finishes the proof of the theorem.

 \end{pf}

 Theorem \ref{Lip is dense in N1pq when q le p} yields the following result.
 \begin{Proposition} \label{Cap pq is outer when q le p plus doubling meas}
 Let $1 \le q \le p<\infty.$ Suppose that $(X,d,\mu)$ satisfies the hypotheses of
 Theorem \ref{Lip is dense in N1pq when q le p}. Then ${\mathrm{Cap}}_{p,q}$ is an outer capacity.

 \end{Proposition}

 In order to prove Proposition \ref{Cap pq is outer when q le p plus doubling meas} we need to state and prove the following proposition, thus generalizing Proposition 1.4 from Bj\"{o}rn-Bj\"{o}rn-Shanmugalingam \cite{BBS}.

 \begin{Proposition} \label{Cap pq is outer for zero sets} (See \cite[Proposition 1.4]{BBS}) Let $1 <p<\infty$ and $1 \le q <\infty.$ Suppose that $(X,d,\mu)$ is a proper metric measure space. Let $E \subset X$ be such that ${\mathrm{Cap}}_{p,q}(E)=0.$ Then for every $\varepsilon>0$ there exists an open set $U \supset E$ with ${\mathrm{Cap}}_{p,q}(U)<\varepsilon.$

 \end{Proposition}

 \begin{pf} We adjust the proof of Proposition 1.4 in Bj\"{o}rn-Bj\"{o}rn-Shanmugalingam \cite{BBS} to the Lorentz setting with some modifications. It is enough to consider the case when $q \neq p.$ Due to the countable subadditivity of ${\mathrm{Cap}}_{p,q}(\cdot)^{1/p}$ we can assume without loss of generality that $E$ is bounded. Moreover, we can also assume that $E$ is Borel. Since ${\mathrm{Cap}}_{p,q}(E)=0,$ we have $\chi_{E} \in N^{1, L^{p,q}}(X,\mu)$ and $||\chi_{E}||_{N^{1,L^{p,q}}}=0.$ Let $\varepsilon \in (0,1)$ be arbitrary. Via Lemma \ref{pq weak upper gradient approximated by upper gradients} and Corollary \ref{Convergence of the quasinorm of the sequence}, there exists $g \in L^{p,q}(X,\mu)$ such that $g$ is an upper gradient for $\chi_{E}$ and $||g||_{L^{p,q}(X,\mu)}<\varepsilon.$ By adapting the proof of the Vitali-Carath\'{e}odory theorem to the Lorentz setting (see Folland
 \cite[Proposition 7.14]{Fol1}) we can find a lower semicontinuous function $\rho \in L^{p,q}(X,\mu)$ such that $\rho \ge g$ and $||\rho-g||_{L^{p,q}(X,\mu)}<\varepsilon.$
 Since ${\mathrm{Cap}}_{p,q}(E)=0,$ it follows immediately that $\mu(E)=0.$ By using the outer regularity of the measure $\mu$ and the absolute continuity of the $(p,q)$-norm, there exists a bounded open set $V \supset E$ such that
 $$||\chi_V||_{L^{p,q}(X,\mu)}+||(\rho+1) \chi_{V}||_{L^{p,q}(X,\mu)}< \frac{\varepsilon}{2}.$$
 Let $$u(x)=\min \bigg \lbrace 1, \inf_{\gamma} \int_{\gamma} (\rho+1) \bigg \rbrace,$$
 where the infimum is taken over all the rectifiable (including constant) curves connecting $x$ to the closed set $X \setminus V.$ We notice immediately that $0 \le u \le 1$ on $X$ and $u=0$ on $X \setminus V.$ By Bj\"{o}rn-Bj\"{o}rn-Shanmugalingam \cite[Lemma 3.3]{BBS} it follows that $u$ is lower semicontinuous on $X$ and thus the set $U=\{ x \in X: u(x)>\frac{1}{2} \}$ is open. We notice that for $x \in E$ and every curve connecting $x$ to some $y \in X \setminus V,$ we have
 $$\int_{\gamma} (\rho+1) \ge \int_{\gamma} \rho \ge \chi_{E}(x)-\chi_{E}(y)=1.$$
 Thus $u=1$ on $E$ and $E \subset U \subset V.$ From \cite[Lemmas 3.1 and 3.2]{BBS}
 it follows that $(\rho+1) \chi_{V}$ is an upper gradient of $u.$ Since $0 \le u \le \chi_{V}$ and $u$ is lower semicontinuous, it follows that $u \in N^{1, L^{p,q}}(X,\mu).$ Moreover, $2u \in {\mathcal{A}}(U)$ and thus
 \begin{eqnarray*}
 {\mathrm{Cap}}_{p,q}(U)^{1/p} &\le& 2 ||u||_{N^{1, L^{p,q}}} \le 2 (||u||_{L^{p,q}(X,\mu)}+||(\rho+1) \chi_{V}||_{L^{p,q}(X,\mu)}) \\
 &\le& 2 (||\chi_V||_{L^{p,q}(X,\mu)} +||(\rho+1) \chi_{V}||_{L^{p,q}(X,\mu)}) < \varepsilon.
 \end{eqnarray*}
 This finishes the proof of Proposition \ref{Cap pq is outer for zero sets}.

 \end{pf}

 Now we prove Proposition \ref{Cap pq is outer when q le p plus doubling meas}.

 \begin{pf} We start the proof of Proposition \ref{Cap pq is outer when q le p plus doubling meas} by showing that every function $u$ in $N^{1, L^{p,q}}(X,\mu)$ is continuous outside open sets of arbitrarily small $p,q$-capacity. (Such a function is called $p,q$-quasicontinuous.) Indeed, let $u$ be a function in $N^{1, L^{p,q}}(X,\mu).$ From Theorem
 \ref{Lip is dense in N1pq when q le p} there exists a sequence $\{u_{j}\}_{j=1}^{\infty}$ of Lipschitz functions on $X$ such that $$||u_{j}-u||_{N^{1, L^{p,q}}}<2^{-2j} \mbox { for every integer } j \ge 1.$$
 For every integer $j \ge 1$ let $$E_j=\{ x \in X: |u_{j+1}(x)-u_{j}(x)|>2^{-j}  \}.$$
 Then all the sets $E_{j}$ are open because the all functions $u_{j}$ are Lipschitz.
 By letting $F=\cap_{j=1}^{\infty} \cup_{k=j}^{\infty} E_k$ and applying the argument from Theorem \ref{Newtonian pq space is Banach} to the sequence $\{u_k \}_{k=1}^{\infty}$ which is Cauchy in $N^{1, L^{p,q}}(X,\mu),$ we see that ${\mathrm{Cap}}_{p,q}(F)=0$ and the sequence $\{u_k\}$
 converges in $N^{1,L^{p,q}}(X,\mu)$ to a function $\widetilde{u}$ whose restriction on $X \setminus F$ is continuous. Thus $||u-\widetilde{u}||_{N^{1, L^{p,q}}}=0$ and hence if we let $E= \{x \in X: u(x) \neq \widetilde{u}(x) \},$ we have ${\mathrm{Cap}}_{p,q} (E)=0.$ Therefore ${\mathrm{Cap}}_{p,q} (E \cup F)=0$ and hence, via Proposition
 \ref{Cap pq is outer for zero sets} we have that $u=\widetilde{u}$ outside open supersets of $E \cup F$ of arbitrarily small $p,q$-capacity. This shows that $u$ is quasicontinuous.

 Now we fix $E \subset X$ and we show that $${\mathrm{Cap}}_{p,q}(E)=\inf \{ {\mathrm{Cap}}_{p,q}(U), E \subset U \subset X, \, U \mbox{ open}  \}.$$ For a fixed $\varepsilon \in (0,1)$ we choose $u \in {\mathcal{A}}(E)$ such that $0 \le u \le 1$ on $X$ and such that
 $$||u||_{N^{1, L^{p,q}}} < {\mathrm{Cap}}_{p,q}(E)^{1/p}+ \varepsilon.$$
 We have that $u$ is $p,q$-quasicontinuous and hence there is an open set $V$ such that
 ${\mathrm{Cap}}_{p,q}(V)^{1/p}<\varepsilon$ and such that $u|_{X \setminus V}$ is continuous. Thus there exists an open set $U$ such that $U \setminus V=\{ x \in X:u(x)>1-\varepsilon \} \setminus V \supset E \setminus V.$
 We see that $U \cup V=(U \setminus V) \cup V$ is an open set containing $E \cup V=(E \setminus V) \cup V.$ Therefore
 \begin{eqnarray*}
 {\mathrm{Cap}}_{p,q}(E)^{1/p} &\le& {\mathrm{Cap}}_{p,q}(U \cup V)^{1/p} \le
 {\mathrm{Cap}}_{p,q}(U \setminus V)^{1/p}+{\mathrm{Cap}}_{p,q}(V)^{1/p} \\
 &\le& \frac{1}{1-\varepsilon} ||u||_{N^{1, L^{p,q}}}+ {\mathrm{Cap}}_{p,q}(V)^{1/p} \\
 &\le& \frac{1}{1-\varepsilon} ({\mathrm{Cap}}_{p,q}(E)^{1/p}+ \varepsilon)+ \varepsilon.
 \end{eqnarray*}
 Letting $\varepsilon \rightarrow 0$ finishes the proof of Proposition \ref{Cap pq is outer when q le p plus doubling meas}.

 \end{pf}

 Theorems \ref{p,q Cap Thm} and \ref{Lip is dense in N1pq when q le p} together with Proposition \ref{Cap pq is outer when q le p plus doubling meas} and Remark \ref{Remarks about the pq capacity} yield immediately the following capacitability result. (See also Appendix II in Doob \cite{Doo}.)

 \begin{Theorem} \label{Cap pq is Choquet when q le p plus doubling meas}
 Let $1<q \le p<\infty.$ Suppose that $(X,d,\mu)$ satisfies the hypotheses of Theorem \ref{Lip is dense in N1pq when q le p}. The set function $E \mapsto {\mathrm{Cap}}_{p,q}(E)$ is a Choquet capacity. In particular, all Borel subsets (in fact, all analytic subsets) $E$ of $X$ are capacitable, that is
 $${\mathrm{Cap}}_{p,q}(E)= \sup \{ {\mathrm{Cap}}_{p,q}(K): K \subset E, \, K \mbox{ \rm{compact}} \}$$
 whenever $E \subset X$ is Borel (or analytic).
 \end{Theorem}

 \begin{Remark} \label{Lip is not dense in N1pinfinity} Counterexample \ref{counterexample conv truncation for N1pinfinity} gives also a counterexample to the density result for $N^{1, L^{p,\infty}}$ in the Euclidean setting for $1<p \le n$ and $q=\infty.$

 \end{Remark}

\end{document}